\documentclass[preprint, times, numbers, sort&compress]{elsarticle}

\usepackage{amssymb}
\usepackage{amsthm}
\usepackage{graphicx}%
\usepackage{multirow}%
\usepackage{amsmath,amssymb,amsfonts}%
\usepackage{mathrsfs}%
\usepackage{dsfont}%
\usepackage[title]{appendix}%
\usepackage{xcolor}%
\usepackage{textcomp}%
\usepackage{manyfoot}%
\usepackage{booktabs}%
\usepackage{algorithm}%
\usepackage{algorithmicx}%
\usepackage{algpseudocode}%
\usepackage{listings}%
\usepackage{url}
\usepackage{makecell}

\newcommand{\vd}{\mathrm{d}}

\theoremstyle{definition}
\newtheorem{teo}{Theorem}

\newtheorem{lem}[teo]{Lemma}

\newtheorem{cor}[teo]{Corollary}
\newtheorem{rem}[teo]{Remark}

\journal{Nuclear Physics B}

\date{}
\begin{document}

\begin{frontmatter}
 \author[first]{Ubong Sam IDIONG\corref{cor1}}
 \affiliation[first]{organization={Department of Mathematics}, address line={Adeyemi Federal University of Education},
                              postcode={351103},
            state ={Ondo},
            country ={NIGERIA}}
            \ead{idiongus@afued.edu.ng}

\author[second]{Unanaowo Nyong BASSEY}
\affiliation[second]{organization={Department of Mathematics},
 address line={University of Ibadan, Ibadan,  Nigeria},
           city={Ibadan},
       postcode={200005},
          state={Oyo State},
          country={NIGERIA}}
\ead{unbassey@yahoo.com}
\author[third]{Olawale Sunday OBABIYI}
\affiliation[third]{organization={Department of Mathematics},
 address line={University of Ibadan, Ibadan,  Nigeria},
           city={Ibadan},
       postcode={200005},
          state={Oyo State},
          country={NIGERIA}}
\ead{so.obabiyi@mail.ui.edu.ng}

\title{ON RADIAL DISTRIBUTION AND QUASI-EXACT SOLVABILITY OF BRIOSCHI-HALPHEN EQUATION}

\begin{abstract}
\noindent The Brioschi-Halphen equation (BHE) is a second order complex differential equation obtained by a two step transformation of the Lam\'e equation. The Lam\'e equation is an equation in Astronomical physics used in the study of motion of planetary bodies. In this paper, the radial part of the BHE for sufficiently large $r$ and the argument limit $2\pi$ is obtained. The asymptotic radial wave function associated with BHE is obtained in terms of canonical polynomials $\mathscr{P}_{n+1},$ and spherical function in $L^{2}(G,\vd\mu), G=SL(2,\mathds{R})$ using point canonical transformation and distributional solution in $\mathscr{C}_{c}^{\infty}(\Omega)$ using Fourier transform method are obtained.
\end{abstract}
\end{frontmatter}
\section{Introduction}
The Brioschi-Halphen equation (BHE) is a second order linear ordinary differential in the complex domain (Poole~\cite{POO}, Chapter IX, \S37, p.163) defined by
\begin{equation}\label{bh}
 \left(4\prod_{s=1}^{3}(w-e_s)D^2+(1-2n)(6w^2-\frac{1}{2}g_2)D+4n(2n-1)w-4B\right)\Psi=0,
\end{equation}
where $D=\frac{\vd}{\vd w},$ $e_{s}=\wp(\omega_{s})$ are the singularities of the Weierstrass polynomial $(\wp'(w))^{2}=4\wp(w)^{3}-g_{2}\wp(w)-g_{3}$, $\omega_{s}\;(s=1,2,3)$ being the half-periods of the Weierstrass elliptic $\wp-$function,  $B$ is the accessory parameter of Lam\'e equation and $n$ is an integer. Lie algebraisation and polynomial solution of BHE have been addressed (see \cite{BI1,BI2}). Distributional solutions of equations with not more than three regular singularities have been addressed using Laplace transform approach (see \cite{KANA,LLL}). In this paper, distributional solutions for radial differential equation with four regular singularities obtained through asymptotic variable separation, are considered using the Fourier transform technique.

The outline of the paper is given in what follows. Section~\ref{FOP} presents formulation of the problem. Section~\ref{mf} presents the mathematical formalism required in obtaining algebraisation of the radial part of  BHE. Section~\ref{QESS}, presents quasi-exact-solvability of the radial BHE. In section~\ref{ESS}, we examine the exact solvability of the radial part of the BHE. In section~\ref{DSS}, the distributional solution of the radial BHE is obtained using Fourier transform approach. The conclusion is drawn for the results obtained in this paper in the final section.
\section{Formulation of Problem}\label{FOP}
In this paper, the method of asymptotic separation of variables (see Estrada and Kanwal~\cite{ESK}; Erd\'{e}lyi~\cite{ERD}) is used to obtain (derive) the radial part of the BHE. To this end, let $w=r\zeta$ with $r\in \mathbb{R}\cup\{\infty\}$ and $\zeta=e^{i\theta},$ where $\theta\in\textrm{[}0,2\pi\textrm{]}, i=\sqrt{-1}$. Let $\Gamma\subset_{open}\mathbb{C},$ then the differential $\vd\Psi$ of a complex valued function $\Psi:\Gamma\hookrightarrow\mathbb{C}$ can be expressed as
\begin{equation}\label{hr1}
  \vd\Psi=\frac{\partial\Psi}{\partial w}\vd w+\frac{\partial\Psi}{\partial \overline{w}}\vd \overline{w},
\end{equation}
where $\overline{w}=r\overline{\zeta}=r\zeta^{-1}$ is the complex conjugate of $w$ [see Hunger~\cite{Hung}, Theorem 3.0.1, p.9]. It is easily seen that
\begin{eqnarray}
  \frac{\partial}{\partial w}&=&\frac{\partial r}{\partial w}\cdot\frac{\partial}{\partial r}+\frac{\partial\zeta}{\partial w}\cdot\frac{\partial}{\partial \zeta}\label{hr2}\\
 \frac{\partial}{\partial \overline{w}}&=&\frac{\partial r}{\partial \overline{w}}\cdot\frac{\partial}{\partial r}+\frac{\partial\zeta}{\partial \overline{w}}\cdot \frac{\partial}{\partial \zeta}.\label{hr3}
\end{eqnarray}
By putting \eqref{hr2} and \eqref{hr3} into \eqref{hr1}, one gets
\begin{eqnarray}
  \vd\Psi &=& \left[\frac{\partial r}{\partial w}\cdot\frac{\partial}{\partial r}+\frac{\partial\zeta}{\partial w}\cdot\frac{\partial}{\partial \zeta}\right]\Psi\vd w+\left[\frac{\partial r}{\partial \overline{w}}\cdot\frac{\partial}{\partial r}+\frac{\partial\zeta}{\partial \overline{w}}\cdot \frac{\partial}{\partial \zeta}\right]\Psi \vd \overline{w} \nonumber\\
   &=& \left[\zeta^{-1}\frac{\partial}{\partial r}+r^{-1}\frac{\partial}{\partial r}\right]\Psi\vd w+\left[\zeta\frac{\partial}{\partial r}-r^{-1}\zeta^2\frac{\partial}{\partial\zeta}\right]\Psi\vd\overline{w}\label{hr4}\\
 \frac{\vd\Psi}{\vd w}  &=&\left[\zeta^{-1}\frac{\partial}{\partial r}+r^{-1}\frac{\partial}{\partial r}\right]\Psi+\left[\zeta\frac{\partial}{\partial r}-r^{-1}\zeta^2\frac{\partial}{\partial\zeta}\right]\Psi\frac{\vd\overline{w}}{\vd w}.\label{hr5}
\end{eqnarray}
Since $$\frac{\vd\overline{w}}{\vd w}=\frac{\frac{\vd\overline{w}}{\vd r}}{\frac{\vd w}{\vd r}}=\frac{\frac{\vd\overline{w}}{\vd \zeta}}{\frac{\vd w}{\vd \zeta}}.$$
It shows that $\frac{\vd\overline{w}}{\vd w}$  has two values namely $\zeta^{-2}$ and $-\zeta^{-2}$ respectively. By substituting these values into equations~\eqref{hr4} and \eqref{hr5} respectively, one obtains
\begin{equation}\label{rader}
  D= 2\zeta^{-1}\frac{\partial}{\partial r}=2r^{-1}\frac{\partial}{\partial \zeta}
\end{equation}
and
\begin{equation}\label{rader2}
  D^{2}=4\zeta^{-2}\frac{\partial^{2}}{\partial r^{2}}=4r^{-2}\frac{\partial^{2}}{\partial\zeta^2}.
\end{equation}
By substituting $w=r\zeta,$ the expressions~\eqref{rader} and \eqref{rader2} into equation~\eqref{bh} one gets
\begin{multline}\label{bhr}
 \bigg[4\prod_{s=1}^{3}(r\zeta-e_s)\left(4\zeta^{-2}\frac{\partial^{2}}{\partial r^{2}} \right) +(1-2n)(6r^{2}\zeta^{2}-\frac{1}{2}g_2)\left(2\zeta^{-1}\frac{\partial}{\partial r}\right).\\
 +4n(2n-1)r\zeta-4B\bigg]\Psi=0.
\end{multline}
Following Estrada and Kanwal~\cite{ESK} with a slight modification, let the function $\Psi$ admit the asymptotic separation of variables $$\Psi(r\zeta) = R(r)\Theta(\zeta)+ o(r^{-1}),\;\;\; \;\; r\rightarrow\infty.$$
Let $\theta\longrightarrow 2\pi,\;\; \zeta\longrightarrow 1$ so that $\Psi \thicksim \lambda R(r),$ where  $ \lambda=\Theta(1)\in \mathds{R}\setminus\{0\}$ and divide \eqref{bhr} through by $4\lambda$. Thus, with these limiting conditions, equation~\eqref{bhr} becomes
 \begin{equation*}
\left(4\prod_{s=1}^{3}(r-e_s)\frac{\vd^{2}}{\vd r^{2}} +\frac{1}{2}(1-2n)(6r^{2}-\frac{1}{2}g_2)\frac{\vd}{\vd r}+n(2n-1)r-B\right)R(r)=0,
 \end{equation*}
which may be rewritten as
 \begin{equation}\label{rdq}
\left(\left(4r^{3}- g_{2}r- g_{3}\right)\frac{\vd^{2}}{\vd r^{2}} -(n-\frac{1}{2})(6r^{2}-\frac{1}{2}g_2)\frac{\vd}{\vd r}+n(2n-1)r-B\right)R(r)=0.
 \end{equation}
 In section~\ref{mf} that follows, we present the formalism for algebraization of a second order differential operator in one real variable which we adopt to write the operator of equation~\eqref{rdq} as an element of the universal enveloping algebra of the Lie algebra $sl(2,\mathbb{R})$. The solvability of the Lie algebraic operator which will allow us to obtain the function $R(r)$ asymptotically in terms of an orthogonal polynomial is discussed in section~\ref{mf}.
\section{Mathematical Formalism}\label{mf}
Following Shifman (\cite{SHIF}, p. 2912),  in order to transform equation \eqref{rdq} into a Lie algebraic equation, we consider a second order differential operator
\begin{equation}\label{op}
  \mathcal{H}:=-\,\frac{1}{2}P_{4}(r)\frac{\vd^2}{\vd r^2}+P_3(r)\frac{\vd}{\vd r}+P_2(r),
\end{equation}
where each $P_\ell(r) \;\;\;(\ell=2,3, 4)$ are polynomials of degree $\ell$ given by
\begin{eqnarray}\label{pc}
  P_4(r) &=& -2[c_{++}r^4-2c_{+0}r^3+(c_{00}-2c_{+-})r^2+2c_{-0}r+c_{--}], \nonumber\\
  P_3(r)&=& -2kc_{++}r^3+(3kc_{+0}-c_{+})r^2 +(-k(c_{00}-2c_{+-})+c_0)r\nonumber\\
  &&\hspace{7cm}+(-kc_{-0}+c_{-}),\nonumber\\
  P_2(r)&=&2jkc_{++}r^2+(-2jkc_{+0}+2jc_{+})r+c_{\ast}.
\end{eqnarray}
Here also, following Shifman (\cite{SHIF}, p. 2913), for the operator $\mathcal{H}$ given in \eqref{op} to be exactly solvable, we assume that $c_{+-}=0$ . Let $\mathfrak{g}$ be the Lie algebra of a certain Lie group $G$. For $1\leq a,b\leq n$ where $n=\dim(\mathfrak{g})$, let $\{\mathcal{J}^{\alpha},\alpha=a,b\}$ denote the generators of the Lie algebra $\mathfrak{g}$ . By rewriting the differential operator in \eqref{op} as the quadratic combination of the generators $\mathcal{J}_{\alpha}$, we get a Lie algebraic operator (cf: Ganguli~\cite{GAN}, p.4)
\begin{equation}\label{pol}
-\mathcal{H}=  \sum_{a,b=0,\pm}c_{a\,b}\,\mathcal{J}_{a}\mathcal{J}_{b}+\sum_{a=0,\pm}c_{a}\,\mathcal{J}_{a}+c_{\ast}.
\end{equation}
Here, $c_{\ast}$ is the Casimir eigenvalue, $c_{a\,b}$ and $c_a$  are complex (or real) numbers and $c_{a\,b}$  form the entries
of the metric $[c_{ab}]_{a,b=0,\pm}$ of dimension $N$ . The generators $\mathcal{J}_{\alpha} (\alpha=a,b)$ of $\mathfrak{g}$ leave invariant the space, say $\mathscr{P}_{n+1}$, of polynomials of degree $n+1$ in the variable $r$, that is, $\mathscr{P}_{n+1} =\langle 1,r,\ldots,r^{n}\rangle$. By construction, the differential operator given in \eqref{pol} also preserves $\mathscr{P}_{n+1}$. The operator representation in \eqref{pol} is known in the
mathematics literature to be isomorphic to the matrix representation of the operator \eqref{op} in finite
dimension [see Brihaye and Godart (\cite{BG}, pp. 5284-5285) and Shifman (\cite{SHIF}, Appendix A, p. 2946)]. It is  called an algebraization of the operator  \eqref{op} and it is an element of the universal enveloping algebra of the Lie algebra $\mathfrak{g}$ of the Lie group $G$ under consideration.

Following Kynsinski~(Kisynski~\cite{KYN}, \S 2, p. 238), the adjoint $\mathcal{H}^{\dag}$ of the Hamiltonian $\mathcal{H}$ in \eqref{pol} is defined as
\begin{equation}\label{adjt}
  \mathcal{H}^{\dag}= - \sum_{a,b=0,\pm}c_{a\,b}\,\mathcal{J}_{a}\mathcal{J}_{b}+\sum_{a=0,\pm}\overline{c}_{a}\,\mathcal{J}_{a}+\overline{c_{\ast}},
\end{equation}
 where $\overline{c}_{a}$ is the complex conjugate of $c_{a}$ in \eqref{pol}.

In this section, operators in one dimension (particularly on real line) is of importance. Take, a Lie group $G=SL(2,\mathbb{R}).$ Then, the Lie algebra $\mathfrak{g}=Lie(G)=sl(2,\mathbb{R}),$  of real matrices with standard basis $H, X^{\pm}$ of order $2\times 2$ and zero traces given by
\begin{equation}\label{gns}
  H=\left(
      \begin{array}{cc}
        \frac{1}{2} & 0 \\
        0 & -\frac{1}{2} \\
      \end{array}
    \right),;\;\;X^{+}=\left(
                        \begin{array}{cc}
                          0 & 1 \\
                          0 & 0 \\
                        \end{array}
                      \right),;\;\;X^{-}=\left(
                                          \begin{array}{cc}
                                            0 & 0 \\
                                            1 & 0 \\
                                          \end{array}
                                        \right)
\end{equation}
satisfying the commutation relation
\begin{equation}\label{gns1}
  [X^{+},X^{-}]=2H, [H,X^{\pm}]=\pm X^{\pm}.
\end{equation}
   Let $\mathfrak{g}_{n}, n\in\mathbb{N}_{0}$ be Lie algebra spanned by first order differential operators depending on non-negative integers $n=2j$
\begin{equation}\label{gno}
\mathcal{J}_{-}=\frac{\vd}{\vd r},\;\;\mathcal{J}_{0}=r\frac{\vd}{\vd r}-j,\;\;\mathcal{J}_{+}=r^2\frac{\vd}{\vd r}-2jr
\end{equation}
which satisfies the same commutation relation as \eqref{gns1}, that is
\begin{equation}\label{gns2}
  [\mathcal{J}_{+},\mathcal{J}_{-}]=2\mathcal{J}_{0},\;\;[\mathcal{J}_{0},\mathcal{J}_{\pm}]=\pm \mathcal{J}_{\pm}.
\end{equation}
Here, $j$ is the spin. Then $\mathfrak{g}_{n}\simeq \mathfrak{g}$. Let $\widehat{\mathfrak{g}}_{n}$ be a central extension of $\mathfrak{g}_n$ by constant functions. Since $\mathfrak{g}_{n}$ is a Lie algebra then any Lie algebraic equation \eqref{pol} for the Lie algebra $\widehat{\mathfrak{g}}_{n}$ is also automatically a Lie algebraic operator for $\mathfrak{g}_n.$ Therefore, the most general second order quasi-exactly-solvable (QES) Hamiltonian in one dimension can be written in the form
\begin{multline}\label{gns3}
-\mathcal{H}=c_{++}(\mathcal{J}_{+})^2+c_{+0}[\mathcal{J}_{+}\mathcal{J}_{0}+\mathcal{J}_{0}\mathcal{J}_{+}]+c_{00}(\mathcal{J}_{0})^{2}+c_{+-}[\mathcal{J}_{+}\mathcal{J}_{-}+\mathcal{J}_{-}\mathcal{J}_{+}]\\
+c_{0-}[\mathcal{J}_{0}\mathcal{J}_{-}+\mathcal{J}_{-}\mathcal{J}_{0}]+c_{--}(\mathcal{J}_{-})^{2}+c_{+}\mathcal{J}_{+}+c_{0}\mathcal{J}_{0}+c_{-}\mathcal{J}_{-}+c_{\ast}.
\end{multline}
The result given in Theorem~\ref{GOL} below was stated in \cite{GKO} without proof. The Theorem is hereby stated with proof for scholastic purpose.
\teo[Gonzalo-Lopez, \emph{et. al.}~\cite{GKO}]\label{GOL} The explicit form of the general second order quasi-exactly-solvable (QES) Hamiltonian \eqref{gns3} is
\begin{multline}\label{cnnf}
-\mathcal{H}=P_{4}(r)\frac{\vd^2}{\vd r^2}+\bigg\{P_{2}(r)-(\frac{2j-1}{2})P_{4}'(r)\bigg\}\frac{\vd}{\vd r}+\bigg\{P_{0}(r)-jP_{2}'(r)+\frac{j(2j-1)}{6}P_{4}''(r)\bigg\}
\end{multline}
where
\begin{eqnarray*}
  P_{4}(r) &=& c_{++}r^{4}+2c_{+0}r^{3}+[c_{00}+2c_{+-}]r^{2}+2c_{0-}r+c_{--}  \\
  P_2(r) &=& c_{+}r^2+c_{0}r+c_{-}\\
  P_{0}(r) &=&\frac{j(j+1)}{3}[c_{00}-4c_{+-}]+c_{*}.
\end{eqnarray*}
\begin{proof}
 Let us assume symmetry, that is, $c_{ab}=c_{ba}$ where $a, b,=0,\pm$. Each term of  $-\mathcal{H}$ are determined as follows.
\begin{eqnarray}
c_{++}\mathcal{J}_{+}^{2}\psi &=&c_{++}\left(r^2\frac{\vd}{\vd r}-2jr\right)\left(r^2\frac{\vd}{\vd r}-2jr\right)\psi\nonumber\\
&=&\left(c_{++}r^{4}\frac{\vd^{2}}{\vd r^{2}}+2(1-2j)c_{++}r^{3}\frac{\vd}{\vd r}-2j(1-2j)c_{++}r^{2}\right)\psi,\label{ag1}\\
c_{00}\mathcal{J}_{0}^{2}\psi &=&c_{00}\left(r\frac{\vd}{\vd r}-j\right)\left(r\frac{\vd}{\vd r}-j\right)\psi\nonumber\\
&=&\left(c_{00}r^2\frac{\vd^2}{\vd r^2}+(2j-1)c_{00}r\frac{\vd}{\vd r} +c_{00}j^2\right)\psi,\label{ag2}\\
c_{--}\mathcal{J}_{-}^{2}\psi&=&c_{--}\frac{\vd^2}{\vd r^2}\psi,\label{ag3}\\
c_{+0}\mathcal{J}_{+}\mathcal{J}_{0}\psi &=& c_{+0} \left(r^2\frac{\vd}{\vd r}-2jr\right) \left(r\frac{\vd}{\vd r}-j\right)\psi\nonumber\\
&=&\left(c_{+0}r^{3}\frac{\vd^2}{\vd r^2}+(1-3j)c_{+0}r^2\frac{\vd}{\vd r}+2j^2c_{+0} r\right),\label{ag4}\\
c_{+0}\mathcal{J}_{0}\mathcal{J}_{+}\psi&=& c_{+0}\left(r\frac{\vd}{\vd r}-j\right)\left(r^2\frac{\vd}{\vd r}-2jr\right)\psi\nonumber\\
&=& \left(c_{+0}r^3\frac{\vd^2}{\vd r^2}+(2-3j)c_{+0}r^2\frac{\vd}{\vd r}+2j(j-1)c_{+0}r\right)\psi,\label{ag5}
\end{eqnarray}
\begin{eqnarray}
c_{0-}\mathcal{J}_{0}\mathcal{J}_{-}\psi &=&c_{0-}\left(r\frac{\vd}{\vd r}-j\right)\frac{\vd\psi}{\vd r}\nonumber\\
&=&\left(c_{0-}r\frac{\vd^2}{\vd r^2}-jc_{0-}\frac{\vd}{\vd r}\right)\psi,\label{ag6}\\
c_{0-}\mathcal{J}_{-}\mathcal{J}_{0}\psi &=& c_{0-}\frac{\vd}{\vd r}\left(r\frac{\vd}{\vd r}-j\right)\psi\nonumber\\
&=&\left(c_{0-}r\frac{\vd^2}{\vd r^2}+(1-j)c_{0-}\frac{\vd}{\vd r}\right)\psi,\label{ag7}\\
c_{+-}\mathcal{J}_{+}\mathcal{J}_{-}\psi&=&c_{+-}\left(r^2\frac{\vd}{\vd r}-2jr\right)\frac{\vd}{\vd r}\nonumber\\
&=& \left(c_{+-}r^2\frac{\vd^2}{\vd r^2}-2jc_{+-}r\frac{\vd}{\vd r}\right)\psi,\label{ag8}\\
c_{+-}\mathcal{J}_{-}\mathcal{J}_{+}\psi &=& c_{+-}\frac{\vd}{\vd r}\left(r^2\frac{\vd}{\vd r}-2jr\right)\psi\nonumber\\
&=&\left(c_{+-}r^2\frac{\vd^2}{\vd r^2}+2(1-j)c_{+-}r\frac{\vd}{\vd r}-2jc_{+-}\right)\psi.\label{ag9}
\end{eqnarray}
\begin{eqnarray}
c_{+}\mathcal{J}_{+}\psi&=&\left(c_{+}r^{2}\frac{\vd}{\vd r}-2jc_{+}r\right)\psi; \label{ag9b}\\
c_{0}\mathcal{J}_{0}\psi&=&\left(c_{0}r\frac{\vd}{\vd r}-c_{0}j\right)\psi;\label{ag9c}\\
 c_{-}\mathcal{J}_{-}\psi&=&c_{-}\frac{\vd\psi}{\vd r}.\label{ag9d}
\end{eqnarray}
Adding the terms \eqref{ag1}-\eqref{ag9d}, one gets
\begin{multline}\label{ag10}
  -\mathcal{H}\psi=\bigg[c_{++}r^{4}+2c_{+0}r^{3}+[c_{00}+2c_{+-}]r^{2}+2c_{0-}r+c_{--}\bigg]\frac{\vd^2\psi}{\vd r^2} \\
  +\bigg[(2j-1)[2c_{++}r^3+3c_{+0}r^2+(2c_{+-}+c_{00})r+c_{0-}]\\
  +c_{+}r^2+c_{0}r+c_{-}\bigg]\frac{\vd\psi}{\vd r} +\big[2j(2j-1)c_{++}r^2\\
  +2j(2j-1)c_{+0}r+c_{00}j^2-2jc_{+-}-j[2c_{+}r+c_{0}]+c_{*}\big]\psi.
\end{multline}
Let the coefficient of $\dfrac{\vd^2}{\vd r^2}$ be denoted by $P_{4}(r)$, that is
\begin{equation}\label{ag10a}
  P_{4}(r)=c_{++}r^{4}+2c_{+0}r^{3}+[c_{00}+2c_{+-}]r^{2}+2c_{0-}r+c_{--}.
\end{equation}
By taking the first and second derivatives of $P_{4}(r)$, one gets
\begin{eqnarray}
  \frac{P_4'(r)}{2} &=&[2c_{++}r^3+3c_{+0}r^2+(2c_{+-}+c_{00})r+c_{0-}] \label{ag10b}\\
 \frac{ P_4''(r)}{6} &=& [2c_{++}r^2+2c_{+0}r+\frac{1}{3}(c_{00}+2c_{+-})].\label{ag10c}
\end{eqnarray}
Let $P_{2}(r)=c_{+}r^2+c_{0}r+c_{-}$ and $P_{2}'(r)=2c_{+}r+c_{0}$. By using equation~\eqref{ag10c}, $P_{2}'(r)$ and the coefficient of $\psi(r)$ in equation~\eqref{ag10} one determines $P_{0}(r)$ as follows
\begin{eqnarray*}
  2j(2j-1)c_{++}r^2 &+&2j(2j-1)c_{+0}r+c_{00}j^2-2jc_{+-}-j[2c_{+}r+c_{0}]+c_{*}\\
  &=&j(2j-1)[2c_{++}r^2+2c_{+0}r]+c_{00}j^2-2jc_{+-}-j[2c_{+}r+c_{0}]+c_{*}\\
  &=& j(2j-1)\bigg[\frac{ P_4''(r)}{6}-\frac{1}{3}(c_{00}+2c_{+-})]\bigg] +c_{00}j^2-2jc_{+-}-jP_{2}'(r)+c_{*}\\
 & =&\frac{j(2j-1)}{6}P_4''(r)-jP_{2}'(r)+\bigg[j^2-\frac{j(2j-1)}{3}\bigg]c_{00}+\frac{2j(2j-1)}{3}c_{+-}+c_{*}\\
  &=&\frac{j(2j-1)}{6}P_4''(r)-jP_{2}'(r)+\frac{j(j+1)}{3}c_{00}+\frac{2j(2j-1)}{3}c_{+-}+c_{*}.
\end{eqnarray*}
The constant part of the last expression is
\begin{equation}\label{ag10d}
  P_{0}(r)=\frac{j(j+1)}{3}c_{00}+\frac{2j(2j-1)}{3}c_{+-}+c_{*}.
\end{equation}
By substituting equations \eqref{ag10a}-\eqref{ag10d} into \eqref{ag10} one obtains
\begin{multline*}
-\mathcal{H}=P_{4}(r)\frac{\vd^2}{\vd r^2}+\bigg\{P_{2}(r)-(\frac{2j-1}{2})P_{4}'(r)\bigg\}\frac{\vd}{\vd r}+\bigg\{P_{0}(r)-jP_{2}'(r)+\frac{j(2j-1)}{6}P_{4}''(r)\bigg\}.
\end{multline*}
Hence the result.
\end{proof}
In Corollary~\ref{cor3.1}, we apply Theorem~\ref{GOL} to carry out the $sl(2)$-algebraization of the radial operator of BHE.
\cor\label{cor3.1} The Lie algebraic version of the radial operator
$$H=\left(4r^{3}- g_{2}r- g_{3}\right)\frac{\vd^{2}}{\vd r^{2}} -(2j-\frac{1}{2})(6r^{2}-\frac{1}{2}g_2)\frac{\vd}{\vd r}+2j(4j-1)r-B.$$
generated by the equation~\eqref{rdq} is given by the symmetric, positive definite closable operator
\begin{equation*}
  -H_1=(4r^3-g_2r-g_3)\frac{\vd^2}{\vd r^2}+\left(\frac{9}{2}(2j-1)r^2+\frac{g_2}{4}\right)\frac{\vd}{\vd r}+7j(2j-1)r-B.
\end{equation*}
  which is defined on the space of polynomials $\mathscr{P}_{2j+1}$ of degree $2j+1$ in the variable $r$.
\begin{proof}
The first task here is to determine the structure constants $c_{ab}$ which allows one to write $H$ in the form found in equation~\eqref{gns3}. Thus by comparison of the coefficients in $H$ with those found in equation~\eqref{ag10} one obtains
$$\begin{array}{cccc}
  c_{++}=0; & c_{+0}=2; & c_{00}=-2c_{+-}=0; \\
 c_{0-}=-\frac{g_2}{2}; &  c_{--}=-g_{3};& c_{0}=0 ; \\
  c_{+}=3(\frac{1}{2}-j); &c_{-}=-(2j-\frac{1}{2})\frac{g_2}{2}=-(j-\frac{1}{4})g_2. &
\end{array}$$

Following a technique similar to the one found in Turbiner (\cite{TA}) (see also Turbiner~\cite{TA1}), we determine the entries of the structure metric $[c_{ab}]_{a,b=0,\pm}$ and the vectors $(c_a)_{a=0,\pm}$ as follows. The structure metric of equation~\eqref{rdq} is given as
\begin{equation}\label{strmet2}
g=[c_{ab}]_{a,b=0,\pm}=\left(
  \begin{array}{ccc}
    c_{++} & c_{+0} & c_{+-} \\
     c_{0+}& c_{00} & c_{0-} \\
    c_{-+} & c_{-0} & c_{--} \\
  \end{array}
\right)=\left(
          \begin{array}{ccc}
            0 & 2 & 0 \\
            2& 0 & -\frac{g_2}{2} \\
            0 & -\frac{g_{2}}{2} & -g_3 \\
          \end{array}
        \right).
\end{equation}
The norm of $g=[c_{ab}]_{a,b=0,\pm}$ is given by
\begin{equation}\label{dstrmet3}
  \|g\|=\det([c_{ab}]_{a,b=0,\pm})=4g_{3}\geqslant 0
\end{equation}
and the modulus of the metric vector $g^{a}=(c_{+},c_{-},c_{0})$ is expressed as
$$|g^{a}|=\sqrt{\left(3(j-\frac{1}{2})\right)^{2}+\left(-(j-\frac{1}{4})g_2\right)^{2}}=\frac{1}{2}\sqrt{9(2j-1)^2+\frac{(4j-1)^2}{4}g_2^2}>0.$$
Here, $j$ is the spin quantum number which is usually defined by $j=\frac{n}{2}$ where $n$ is an integer provided the quantum Hamiltonian under consideration has discrete spectrum (see Olver~\cite{OLP} , p.109).

The operator $H$ in \eqref{rdq} takes the canonical form
\begin{multline}\label{og1}
  -H_1=2\mathcal{J}_{+}\mathcal{J}_{0}+2\mathcal{J}_{0}\mathcal{J}_{+}-\frac{g_2}{2}\mathcal{J}_{0}\mathcal{J}_{-}-\frac{g_2}{2}\mathcal{J}_{-}\mathcal{J}_{0}-g_3(\mathcal{J}_{-})^2+(\frac{3}{2}-3j)\mathcal{J}_{+}\\
  +(j-\frac{1}{4})g_2\mathcal{J}_{-}-B.
\end{multline}
 The explicit form of the Lie algebraic Hamiltonian using equation~\eqref{cnnf} is
\begin{equation}\label{cfm}
  -H_1=(4r^3-g_2r-g_3)\frac{\vd^2}{\vd r^2}+\left(\frac{9}{2}(2j-1)r^2+\frac{g_2}{4}\right)\frac{\vd}{\vd r}+7j(2j-1)r-B.
\end{equation}
The  symmetric nature of the structure matrix in \eqref{strmet2} and the value of its determinant in \eqref{dstrmet3} show that the operator $-H_1$  is  symmetric elliptic densely defined operator on $L^2(G,\mu)$ since the domain $dom(H_1)=\mathscr{P}_{n+1} $ is dense in $L^2(G,\mu)$ and $\mathscr{P}_{n+1}$ is invariant under $H_1$.
 \end{proof}
Since $j$ and $g_2$ are real numbers, by using \eqref{adjt}, then the adjoint of $-H_1$ is obtained as
\begin{multline}\label{og1a}
  (-H_1)^{\dag}=2\mathcal{J}_{+}\mathcal{J}_{0}+2\mathcal{J}_{0}\mathcal{J}_{+}-\frac{g_2}{2}\mathcal{J}_{0}\mathcal{J}_{-}-\frac{g_2}{2}\mathcal{J}_{-}\mathcal{J}_{0}-g_3(\mathcal{J}_{-})^2-(\frac{3}{2}-3j)\mathcal{J}_{+}\\
  -(j-\frac{1}{4})g_2\mathcal{J}_{-}+B.
\end{multline}
This shows that $-H_1$ is not a self-adjoint map. By a result in Akhiezer and Glazman (\cite{AKG}, p. 80, vol. 1) the adjoint $(-H_1)^{\dag}$ of the operator $-H_1$ is closed since $H_{1}$ is densely defined in the Hilbert space $L^2(G,\vd\mu)$.
\ \\
\section{Quasi-Exact Solvability}\label{QESS}
The operator $H_1$ is a quasi-exactly solvable (cf: Shifman~\cite{SHIF}, p.2908). Recall that a linear differential operator $T$ in a Hilbert space $\mathscr{H}$ is called quasi-exactly solvable (QES) if it leaves  invariant a non-trivial finite-dimensional  subspace $\mathscr{P}\subset\mathscr{H}$. That is, $$T\mathscr{P}\subset\mathscr{P},\;\mathscr{P}=\langle\varphi_0, \varphi_1,\varphi_2,\ldots,\varphi_n\rangle, \varphi_i\in\mathscr{H}.$$
Then, the first $n$ eigenvalues and corresponding eigenfunctions can be obtained exactly by diagonalizing the corresponding matrix of the restricted action of $T$ to the subspace $\mathscr{P}.$  By Turbiner (cf: Turbiner~\cite{TUB}, Lemma 2.2, p.13) a QES operator $T\in U_{s\ell(2,\mathbb{R})}$ has no terms in $\mathcal{J}_{+}$, positive grading, if and only if it is an exactly solvable operator. Let the operator $H_1$ be restricted to $\mathscr{P}_{n+1}$ is obtained. There are three well-known techniques of solving QES operator equations, namely, Bethe ansatz method~\cite{KAR}, constraint polynomial approach \cite{MOR} and canonical polynomial technique \cite{PZB}. In this work, the latter approach which is the most recent and easiest is adopted.  Using canonical polynomial approach,  when the monomial $r^{n}$ is acted upon by the operator $H_{1},$ it generates a recurrence equation which is solved using Jacobi tri-diagonal matrix (see \cite{ELM}). Since $\mathscr{P}_{n+1}$ is invariant under the action of the operator $H_{1}$ one gets
\begin{eqnarray*}
  -H_{1}r^{k} &=& (4r^3-g_2r-g_3)k(k-1)r^{k-2}+ \left(\frac{9}{2}(2j-1)r^2+\frac{g_2}{4}\right)kr^{k-1}\\
  &&+7j(2j-1)r^{k+1}-Br^k\\
  &=&4k(k-1)r^{k+1}-g_2k(k-1)r^{k-1}+\frac{g_2}{4}kr^{k-1}-g_3k(k-1)r^{k-2}+ \frac{9}{2}(2j-1)kr^{k+1}\\
  &&+7j(2j-1)r^{k+1}-Br^k\\
   &=& \left[4k(k-1)+\frac{9}{2}(2j-1)k+7j(2j-1)\right]r^{k+1}-Br^{k}\\
   && -\frac{g_2}{4}k\left[4(k-1)+1\right]r^{k-1}-k(k-1)g_3r^{k-2}.
\end{eqnarray*}
Thus the entries of the tri-diagonal Jacobi matrix associated with $H_{1}$ are
\begin{eqnarray}
  \tau_{k,k+1} &=& 4k(k-1)+\frac{9}{2}(2j-1)k+7j(2j-1)\label{qr1a} \\
  \tau_{k,k}&=& -B\label{qr1b}\\
  \tau_{k,k-1} &=& -\frac{g_2}{4}k\left[4k-3\right]\label{qr1c}\\
 \tau_{k,k-2} &=& -k(k-1)g_3.\label{qr1d}
\end{eqnarray}
 In what follows, the technique of gauge transformation (cf:Gonzalo-Lopez \emph{et. al.}, ~\cite{GKO}) is applied to eigenvalue equation involving $H_{1}$. The statement of the result will serve as a guide.
 \teo[cf:\cite{GKO}, Theorem~3, Statement only]\label{GKOT} Let
$$-\mathit{L}=\mathcal{P}(r)\frac{\vd^{2}}{\vd r^{2}}+\mathcal{Q}(r)\frac{\vd}{\vd r}+\mathcal{R}(r),\;\;r\in\mathbb{PR}^{1}$$
be a second order ordinary differential operator with $\mathcal{P}(r)> 0$ which satisfies the equation
$$L[\psi(r)]=0.$$
Then by variable transformation
\begin{equation*}
w=\varphi(r)=\int_{r}^{\infty}\left(\mathcal{P}(u)\right)^{-\frac{1}{2}}\vd u
\end{equation*}
and the gauge function
\begin{equation*}
\mu(r)=(\mathcal{P}(r))^{-\frac{1}{4}}\exp\left(\frac{1}{2}\int_{}^{r}\frac{\mathcal{Q}(u)}{\mathcal{P}(u)}\vd u\right)
\end{equation*}
transforms $L$ into Schr\"{o}dinger form by the transformation
\begin{equation}\label{qh3}
\mu(r)\cdot L\cdot \mu(r)^{-1}=- \frac{\vd^{2}}{\vd w^{2}}+V(w)=:\mathcal{S}
\end{equation}
where the potential is given by
\begin{equation*}
V(w)=-\frac{3\mathcal{P}'(r)^2-8\mathcal{P}'(r)\mathcal{Q}(r)+4\mathcal{Q}^2(r)}{16\mathcal{P}(r)}-\frac{1}{4}\mathcal{P}''(r)+\frac{1}{2}\mathcal{Q}'(r)-\mathcal{R}(r)
\end{equation*}
and the eigenfunction  $\psi(r)$ of $L$ gives the eigenfunction of the resulting Schr\"{o}dinger operator in expression~\eqref{qh3}  in its coordinate as
\begin{equation*}
\widetilde{\psi}(w)=\mu(\varphi^{-1}(w))\psi(\varphi^{-1}(w)).
\end{equation*}
In what follows, Theorem~\ref{GKOT} is applied to $H_1$ to examine its quasi-exact solvability.
\teo The operator $H_1$ by gauge transformation in Theorem~\ref{GKOT} yields a Schr\"{o}dinger operator
\[\mathcal{S}:=-\frac{\vd}{\vd w}+V(w),\]
where $w(r)=\wp^{-1}(r)$ is the pull-back of the Weierstrass elliptic $p$-function.
and
\[V(w)=\frac{(12r-g_2)[36r(1-(2j-1)r) -5g_2]}{16(4r^3-g_2r-g_3)}+\frac{(18(2j-1)r^2+g_2)^2}{64(4r^3-g_2r-g_3)}-\frac{1}{2}[28j^2+32j-3]r+B.\]
The associated $n^{\mathrm{th}}$ radial wave function $R_{n}(r)$ is given in terms of polynomial of degree $n, \mathcal{P}_{n}$
\[ R_{n}(r)= 2^{-\frac{n}{2}} \prod_{s=1}^{3}(r-e_s)^{\eta_s-\frac{n}{4}}\mathcal{P}_{n}(r)\]
where
\[ \mathcal{P}_{n}(r)=[1+\sum_{m=1}^{n}\mu_{j}^{(m)}r^{m}],\]
 \[\mu_{j}^{(m)}=\frac{\sum_{k=0}^{m}[\tau_{k,0}+\tau_{k,1}\mu_{j}^{(1)}\ldots+\tau_{k,m-1}\mu_{j}^{(m-1)}]}{\sum_{k=0}^{m}\tau_{k,m}}, j=\frac{n}{2}\]
 and $\tau_{k,0},\tau_{k,1},\ldots,\tau_{k,m}$ are entries of the Jacobi tri-diagonal matrix associated with $H_1.$
\begin{proof}
  Let
 \begin{equation}\label{qr4}
   R_{n}(r)=g(r) \mathcal{P}_{n}(r(w)), \;\;\mathrm{where}\;\;\mathcal{P}_n\in \mathscr{P}_{n+1}=\langle 1, r, r^2, \ldots, r^{n}\rangle
 \end{equation}
 and $g(r)$ is the gauge function. In this case, $H_{1}$ is equivalent to the operator $\mathcal{S}$ under the transformation
 \begin{equation}\label{qr5}
   \mathcal{S}= g(r)^{-1}\circ (-H_{1})\circ g(r).
 \end{equation}
 It should be remarked here that the transformed operator $\widetilde{\mathbb{H}}$ preserves the eigenvalues $E_j$
 and leaves the polynomial space $\mathscr{P}_{m+1}$ invariant. Thus, $\widetilde{\mathbb{H}}\in U_{s\ell(2,\mathbb{R})}$ according to Burnside's theorem (see Weyl~\cite{HW} and Panahia \emph{et. al.}~\cite{PZB}).  Following Gonzalo-Lopez \emph{et. al} (\cite{GKO}, Theorem 3.0, p.117), $\mathcal{S}$ is defined by
 \begin{equation}\label{qr6}
   \mathcal{S}:= -\frac{\vd^2}{\vd w^2}+V(w)
 \end{equation}
 where $V(w)$ is the exactly solvable potential given by
 $$V(w)=\, \frac{3(P'_3)^2-8P'_3P_2+4P_2^2}{16P_3}-\frac{1}{4}P_3''+\frac{1}{2}P'_2-P_1,$$
 provided that
 $$g(r)=(P_3(r))^{-\frac{2j}{4}}\exp\left(\frac{1}{2}\int^{r}\frac{P_2(u)}{P_3(u)}\vd u\right),$$
 $P_{3}(r)=4r^3-g_2r-g_3=4(r-e_1)(r-e_2)(r-e_3)> 0$ (this is a possibility when $r> \max (e_1,e_2, e_3),$ since $e_1, e_2, e_3$ are real numbers which are comparable ) and $$w=\varphi(r)=\int_{r}^{\infty}\frac{\vd u}{\sqrt{P_3(u)}}=\int_{\wp(r)}^{\infty}\frac{\vd u}{\sqrt{4u^3-g_2u-g_3}}=\wp^{-1}(r),$$
 where $\wp^{-1}$ is the pull-back of the Weierstrass elliptic $\wp$-function. Recall that in $H_{1}$,
 \begin{eqnarray*}
   P_2(r) &=& \frac{9}{2}(2j-1)r^2+\frac{g_2}{4}, \\
   P_1(r) &=&  7j(2j-1)r-B.
 \end{eqnarray*}
 The gauge function $g(r)$ is thus
 \begin{eqnarray}
   g(r) &=& 4^{-\frac{2j}{4}}\prod_{s=1}^{3} (r-e_s)^{-\frac{j}{2}}\exp\left(\frac{1}{2}\int^{r}\frac{P_2(u)}{P_3(u)}\vd u\right)\nonumber\\
   &=& 2^{-j}\prod_{s=1}^{3}(r-e_s)^{-\frac{j}{2}}\exp\left(\frac{1}{2}\int^{r}\frac{\frac{9}{2}(2j-1)r^2+\frac{g_2}{4}}{4u^3-g_2u-g_3}\vd u\right)\nonumber\\
   &=&2^{-j}\prod_{s=1}^{3}(r-e_s)^{-\frac{j}{2}}\exp\left(\frac{1}{2}\int^{r}\frac{\frac{9}{2}(2j-1)r^2+\frac{g_2}{4}}{4(r-e_1)(r-e_2)(r-e_3)}\vd u\right)\nonumber\\
   &=& 2^{-j} \prod_{s=1}^{3}(r-e_s)^{\eta_s-\frac{j}{2}},\label{qr7}
 \end{eqnarray}
 where
 \begin{eqnarray*}
   \eta_1 &=& \frac{\frac{9}{2}(2j-1)e_{1}^2+\frac{g_2}{4}}{4(e_1-e_2)(e_1-e_3)}, \\
   \eta_2 &=& \frac{\frac{9}{2}(2j-1)e_{2}^2+\frac{g_2}{4}}{4(e_2-e_1)(e_2-e_3)},\;\text{and} \\
   \eta_3 &=& \frac{\frac{9}{2}(2j-1)e_{3}^2+\frac{g_2}{4})}{4(e_3-e_2)(e_3-e_1)}.
 \end{eqnarray*}
 Recall that
 \begin{equation*}
   \begin{array}{ccc}
   P_3(r)=4r^3-g_2 r -g_3,\;\;& P'_3(r)=3r^2-g_2,\;\; & P''_3(r)=6r \\
   P_2(r)= \frac{9}{2}(2j-1)r^2+\frac{g_2}{4},\;\; & P'_2 =9(2j-1)r,\;\;& P_1(r)=7j(2j-1)r-B
 \end{array}
 \end{equation*}
 then the potential $V(w)$ is obtained as
 \begin{eqnarray*}
   V(w) &=&  \frac{3(P'_3)^2-8P'_3P_2+4P_2^2}{16P_3}-\frac{1}{4}P_3''+\frac{1}{2}P'_2-P_1\\
    &=& \frac{(12r-g_2)[36r(1-(2j-1)r) -5g_2]}{16(4r^3-g_2r-g_3)}+\frac{(18(2j-1)r^2+g_2)^2}{64(4r^3-g_2r-g_3)}\\
    &&-\frac{1}{2}[28j^2+32j-3]r+B.
 \end{eqnarray*}
 Since the gauge function has been obtained, the gauge transformation in \eqref{qr5} can be re-written as
 \begin{equation}\label{qr8}
   R_{n}(r)= 2^{-j} \prod_{s=1}^{3}(r-e_s)^{\eta_s-\frac{j}{2}}\mathcal{P}_{n}(r), \;\; \mathrm{where}\;\; \mathcal{P}_{n}(r)=\sum_{m=0}^{n}a_{m}r^{m}
 \end{equation}
 and $n =2j$ with ($n=0, 1,2,\ldots$). The next task is to determine the coefficients $a_{m} \;(m=0, 1,2, \ldots, n)$ using the tri-diagonal Jacobi matrix form in \eqref{qr1a}-\eqref{qr1d}. For convenience, let the table of entries of the tri-diagonal Jacobi matrix be determined for values $m=0,1, 2,3$ (see Table~\ref{tabqr}). It is noteworthy to see that $\tau_{m,m}=0, \forall\; m=0,\ldots,n$. It is now necessary to look at the nature of Jacobi tri-diagonal matrices, their corresponding eigenvalues as well as eigenfunctions for each case.
  \begin{table}[h]
   \centering
   \begin{tabular}{||p{.5cm}||p{7.5cm}||p{2.5cm}||p{2.5cm}||}
   \hline
   \hline
  $m$ & $\tau_{m,m+1}$ & $\tau_{m,m-1}$&$\tau_{m,m-2}$ \\
\hline
\hline
    $0$ & $\tau_{0,1}=7j(2j-1)$  & $\tau_{0,-1}=0$ & $\tau_{0,-2}=0$  \\
\hline
\hline
    $1$& $\tau_{1,2}=(7j+\frac{9}{2})(2j-1)$ & $\tau_{1,0}=-\frac{g_2}{4}$  &$\tau_{1,-1}=0$\\
\hline
\hline
   $2$&  $\tau_{2,3}=8+9(2j-1)+7j(2j-1)$& $\tau_{2,1}=-\frac{5}{2}g_2$ &$\tau_{2,0}=-2 g_3$\\
\hline
\hline
 \end{tabular}
   \caption{Table of values for $\tau_{m,m+1},\tau_{m,m-1},\tau_{m,m-2}, m=0,1,2,\ldots$}\label{tabqr}
 \end{table}
 \ \\
 We also note here that $\tau_{m,m}=-B, (m=0,1,2,\ldots, n)$.
 \begin{description}
   \item[case (i): ($m=0$).] In this case, $\mathcal{S}$ possesses the invariant subspace $\mathscr{P}_1$ spanned by the basis $\{1\}$, thus, the function $\mathcal{P}_0(r)=a_0$. The matrix equation corresponding $a_0$ is given by the $1\times 1$ matrix equation given by
       $$T_{1}A_{1}=[\tau_{0,0}][a_0]=0\implies -Ba_0=0.$$
       Since $a_0\neq 0$ for ground state eigenfunction $B=0$.
    Thus, by the gauge transformation in \eqref{qr8}, we get the ground state eigenfunction
       $$R_{0}(r)=g(r)\mathcal{P}_{0}(r)=a_0g(r)= 2^{-j}a_{0} \prod_{s=1}^{3}(r-e_s)^{\eta_s-\frac{j}{2}}.$$
   \item[case (ii): ($m=1$).] In this case, $\mathcal{S}$ possesses the invariant subspace $\mathscr{P}_1$ spanned by the basis $\{1,r\}$, thus, the function $\mathcal{P}_1(r)=a_0+a_1r$. The matrix equation corresponding $a_0$ is given by the $2\times 2$ matrix equation given by
      \begin{equation}\label{qr9}
        T_{2}A_{2}=\left(
                      \begin{array}{cc}
                        \tau_{0,0} & \tau_{0,1} \\
                        \tau_{1,0} & \tau_{1,1} \\
                      \end{array}
                    \right)\left(
                             \begin{array}{c}
                               a_0 \\
                               a_1 \\
                             \end{array}
                           \right)=0.
      \end{equation}
       To solve for the values $B$, one notes that since $A_2\neq 0$ then the determinant of $T_2$ must be zero. Therefore,
       \begin{equation*}
         \det T_2 = \left|
                        \begin{array}{cc}
                          \tau_{0,0} & \tau_{0,1} \\
                          \tau_{1,0} & \tau_{1,1} \\
                        \end{array}
                      \right| = \tau_{0,0}\tau_{1,1}-\tau_{1,0}\tau_{0,1}=0.
       \end{equation*}
       This yields a non-zero eigenvalue
       \begin{eqnarray*}
         B^2+\frac{g_2}{4}7j(2j-1) &=& 0 \\
         B_{\pm} &=& \pm\frac{\sqrt{7j(1-2j)g_2}}{2}.
       \end{eqnarray*}

       Next, we solve for $a_1$ in equation~\eqref{qr9}.
       Equation~\eqref{qr9} can be expressed as
        \begin{equation}\label{qr9b}
          \tau_{0,0}a_0+\tau_{0,1}a_1=0;\;\;\tau_{1,0}a_0+\tau_{1,1}a_1=0.
        \end{equation}
        By adding up the two equations one gets
        $$(\tau_{0,1}+\tau_{1,1})a_1=-(\tau_{0,0}+\tau_{1,0})a_0\implies a_1=-\;\frac{(\tau_{0,0}+\tau_{1,0})}{\tau_{0,1}+\tau_{1,1}}a_0.$$
        Thus,
        $$a_1=-\left(\frac{g_2+4B}{4[7j(2j-1)-B]}\right)a_0.$$
        Thus, by the gauge transformation in \eqref{qr8}, we get the first state eigenfunction
       \begin{eqnarray*}
         R_{1}(r)&=&g(r)\mathcal{P}_{1}(r)\\
         &=&[a_0+a_1r]g(r)\\
         &=& 2^{-j}a_{0}\left[1-\left(\frac{g_2+4B}{4[7j(2j-1)-B]}\right)r\right] \prod_{s=1}^{3}(r-e_s)^{\eta_s-\frac{j}{2}}.
       \end{eqnarray*}
       Let $\displaystyle\mu_{j}^{(1)}=\frac{g_2+4B}{4[7j(2j-1)-B]}$ then $a_1=\mu_j^{(1)}a_0$
       and
       $$R_{1}(r)=2^{-j}a_0[1+\mu_{j}^{(1)}r]\prod_{s=1}^{3}(r-e_s)^{\eta_s-\frac{j}{2}}.$$
   \item[case (iii): ($m=2$)] In this case, $\mathcal{S}$ possesses the invariant subspace $\mathscr{P}_1$ spanned by the basis $\{1,r, r^2\}$, thus, the function $\mathcal{P}_2(r)=a_0+a_1r+a_2r^2$. The matrix equation corresponding $a_0$ is given by the $3\times 3$ matrix equation given by
      \begin{equation}\label{qr9}
        T_{3}A_{3}=\left(
                      \begin{array}{ccc}
                        \tau_{0,0} & \tau_{0,1} & \tau_{0,2} \\
                        \tau_{1,0} & \tau_{1,1} &\tau_{1,2} \\
                        \tau_{2,0}& \tau_{2,1}&\tau_{2,2}
                      \end{array}
                    \right)\left(
                             \begin{array}{c}
                               a_0 \\
                               a_1 \\
                               a_2
                             \end{array}
                           \right)=0.
      \end{equation}
       To solve for $B$ in this case, $A_3\neq 0$, thus  $$\det[T_3]=0.$$
       This implies that
       \begin{multline}\label{qr10}
         \tau_{0,0}\tau_{1,1}\tau_{2,2}-\tau_{2,1}\tau_{0,0}\tau_{1,2}
         -\tau_{0,1}\tau_{1,0}\tau_{2,2}+\tau_{0,1}\tau_{2,0}\tau_{1,2}\\
         +\tau_{0,2}\tau_{1,0}\tau_{2,1}-\tau_{0,2}\tau_{1,1}\tau_{2,0}=0.
       \end{multline}
       Since $\tau_{0,0}=\tau_{1,1}=\tau_{2,2}=E_j$, equation~\eqref{qr10} becomes
       \begin{multline}\label{qr10b}
         B^{3}-B(\tau_{2,1}\tau_{1,2}+\tau_{0,1}\tau_{1,0}\tau_{0,2}\tau_{2,0})      +\tau_{0,1}\tau_{2,0}\tau_{1,2}+\tau_{0,2}\tau_{1,0}\tau_{2,1}=0.
       \end{multline}
       By standard formula for solving cubic polynomials (cf: Abramowitz and Stegun~\cite{AS}, \S 3.8.2, p.17)
       $z^3+b_2z^2+b_1z+b_0=0$ has the roots
       \begin{eqnarray*}
         z_1 &=& (s_++s_{-})-\frac{b_2}{3}, \\
         z_2 &=& -\frac{1}{2}(s_++s_{-})-\frac{b_2}{3}+i\frac{\sqrt{3}}{2}(s_+-s_{-}),\\
         z_3 &=&  -\frac{1}{2}(s_++s_{-})-\frac{b_2}{3}-i\frac{\sqrt{3}}{2}(s_+-s_{-}),\\
       \end{eqnarray*}
       where,
       \begin{eqnarray*}
         s_\pm &=& [t\pm(q^3+t^2)^{\frac{1}{2}}]^\frac{1}{3}, \\
          q &=& \frac{1}{3}b_1-\frac{1}{9}b_2^2,\\
         t &=& \frac{1}{6}(b_1b_2-b_0)-\frac{1}{27}b_2^3.
       \end{eqnarray*}
       Thus, by setting $$b_2=0, b_1=-(\tau_{2,1}\tau_{1,2}+\tau_{0,1}\tau_{1,0}\tau_{0,2}\tau_{2,0}), b_{0}=\tau_{0,1}\tau_{2,0}\tau_{1,2}+\tau_{0,2}\tau_{1,0}\tau_{2,1}$$
       \begin{eqnarray*}
         q &=& -\frac{1}{3}(\tau_{2,1}\tau_{1,2}+\tau_{0,1}\tau_{1,0}\tau_{0,2}\tau_{2,0}),\\
         t &=& -\frac{1}{6}(\tau_{0,1}\tau_{2,0}\tau_{1,2}+\tau_{0,2}\tau_{1,0}\tau_{2,1}),   \\
         s_\pm &=& \left[-\frac{1}{6}(\tau_{0,1}\tau_{2,0}\tau_{1,2}+\tau_{0,2}\tau_{1,0}\tau_{2,1})\right.\\
         && \left. \pm\left(-\frac{1}{27}(\tau_{2,1}\tau_{1,2}+\tau_{0,1}\tau_{1,0}\tau_{0,2}\tau_{2,0})^3+\frac{(\tau_{0,1}\tau_{2,0}\tau_{1,2}+\tau_{0,2}\tau_{1,0}\tau_{2,1})^2}{36}\right)^{\frac{1}{2}}\right]^{\frac{1}{3}}.
       \end{eqnarray*}
       Therefore the roots of equation~\eqref{qr10b} are
       \begin{eqnarray*}
         B_{0} &=& (s_++s_{-}),\\
        B_{+} &=& -\frac{1}{2}(s_++s_{-})+i\frac{\sqrt{3}}{2}(s_+-s_{-}),\\
         B_{-} &=&  -\frac{1}{2}(s_++s_{-})-i\frac{\sqrt{3}}{2}(s_+-s_{-}).\\
       \end{eqnarray*}
Solving for $a_2$ in equation~\eqref{qr9}, one gets
$$ a_2= \frac{\displaystyle\sum_{m=0}^{2}[\tau_{m,0}+\tau_{m,1}\mu_j^{(1)}]a_0}{\displaystyle\sum_{m=0}^{2}\tau_{m,2}}=\mu_j^{(2)}a_0 $$
where $$\mu_j^{(2)}=-\frac{[4B+g_2+8g_3+[28j(2j-1)-4B-10g_2]\mu_{j}^{(1)}]}{(28j+18)(2j-1)-4B}.$$
 Thus, by the gauge transformation in \eqref{qr8}, we get the second state eigenfunction
 \begin{eqnarray*}
         R_{2}(r)&=&g(r)\mathcal{P}_{2}(r)\\
         &=&a_0[1+\mu_j^{(1)}r+\mu_j^{(2)}r^2]g(r)\\
         &=&2^{-j}a_0[1+\mu_j^{(1)}r+\mu_j^{(2)}r^2] \prod_{s=1}^{3}(r-e_s)^{\eta_s-\frac{j}{2}}.
  \end{eqnarray*}
   \item[case (v): ($m=2j$).] In this case, the operator $\mathcal{S}$ has a finite-dimensional
invariant subspace $\mathscr{P}_ {2j+1}$, which is spanned by the basis $\{r^{m}|m=0,1,2,\ldots,2j\}$. By the finite polynomial $\displaystyle\mathcal{P}_{2j}(r)=\sum_{m=0}^{2j} a_{m}r^{m}$ it is possible to obtain the tri-diagonal Jacobi matrix where $2j=n$ as follows
\begin{equation}\label{qr14}
 T_{n+1}A_{n+1}= \left(
    \begin{array}{ccccccc}
      B & \tau_{0,1} & 0 & 0 &\ldots  &&0 \\
      \tau_{1,0} & B & \tau_{1,2} & 0 & \ldots &&\vdots\\
      \tau_{2,0} & \tau_{2,1} & B & \tau_{2,3} & \ldots &&  \\
      0 & \tau_{3,1} & \tau_{3,2} &B&\ddots&&  \\
      \vdots & 0 & \tau_{4,2} & \tau_{4,3} &\ddots&& 0\\
      0& 0&\ldots &&&&\tau_{n-1,n}\\
      0 & \ldots &&&\ldots&\tau_{n,n-1}&B
    \end{array}
  \right)\left(
           \begin{array}{c}
             a_0 \\
             a_1 \\
             a_2 \\
             \\
             \\
             \vdots \\
             a_{n-1}\\
             a_{n} \\
           \end{array}
         \right)=0
\end{equation}
so that one can determine its eigenvectors $\mathbf{B}=[B^{(n)}_{jl}] (l=0,\ldots,n)$ and the general formula for obtaining the coefficients
$a_m$ can be obtained using recursive relation obtained from the tri-diagonal Jacobi matrix since $H_{1}$ preserves the polynomial $\mathcal{P}_{n}$. Therefore,
\begin{eqnarray*}
- H_{1}\mathcal{P}_n(r)&=&-\sum_{m=0}^{n}a_{m}H_{1}r^{m}\\
  &=&\sum_{m=0}^{n}a_{m}[\tau_{m,m+1}r^{m+1}+\tau_{m,m-1}r^{m-1}+\tau_{m,m-2}r^{m-2}-Br^{m}].\\
  \end{eqnarray*}
Since $a_{m-2}=a_{m-1}=0=\tau_{m,m-1}=\tau_{m,m-2}$  for $m=0,1$ one gets.
\begin{equation*}
  -H_{1}\mathcal{P}_{n}=\sum_{m=0}^{n}\left([\tau_{m-1,m}-B]a_{m-1}+\tau_{m+1,m}a_{m+1}+\tau_{m+2,m}a_{m+2}\right)r^{m}=0.
\end{equation*}
Thus, the 3-term recursive relation obtained is
 $$[\tau_{m-1,m}-B]a_{m-1}+\tau_{m+1,m}a_{m+1}+\tau_{m+2,m}a_{m+2}=0,\;\;m\ge 1,$$
 where
 \begin{eqnarray*}
  \tau_{m-1, m} &=& 4(m-1)(m-2)+\frac{9}{2}(2j-1)(m-1)+7j(2j-1) \\
  \tau_{m+1, m} &=& -\frac{g_2}{4}(m+1)\left[4m+1\right]\\
 \tau_{m+2, m} &=& -(m+2)(m+1)g_3.
  \end{eqnarray*}
  Therefore, the general form of the eigenfunction is
  \begin{equation*}
    R_{n}(r)=2^{-j}a_0[1+\sum_{m=1}^{2j}\mu_{j}^{(m)}r^{m}] \prod_{s=1}^{3}(r-e_s)^{\eta_s-\frac{j}{2}}.
  \end{equation*}
  Here $$\mu_{j}^{(m)}=\frac{\sum_{k=0}^{m}[\tau_{k,0}+\tau_{k,1}\mu_{j}^{(1)}\ldots+\tau_{k,m-1}\mu_{j}^{(m-1)}]}{\sum_{k=0}^{m}\tau_{k,m}}.$$
  \end{description}
  To determine the eigenvalues $B$ we use the determinant formula since $A_{n+1}\neq 0$, we have
  \begin{equation}\label{dform}
    D_{n+1}=\det(T_{n+1})=0,
  \end{equation}
  where
  \begin{equation}\label{dform1}
    D_{n+1}=\tau_{n,n}D_{n}-\tau_{n,n-1}\tau_{n-1,n}D_{n-1}
  \end{equation}
  with initial conditions $D_{-1}=0, \;D_{0}=1,\;D_{1}=\tau_{0,0}=B$.
  Since $\tau_{n,n}=B$ for all $n=2j\in\mathbb{N}_{0},$ we rewrite equation~\eqref{dform1} as
  \begin{equation}\label{dform2}
     D_{n+1}=BD_{n}-\tau_{n,n-1}\tau_{n-1,n}D_{n-1}.
  \end{equation}
  The characteristic equation associated with equation\eqref{dform2} is given by
  $$\lambda^2-B\lambda+\tau_{n,n-1}\tau_{n-1,n}=0$$
  with roots
  $$\lambda_{\pm}=\frac{B\pm\sqrt{B^2-4\tau_{n,n-1}\tau_{n-1,n}}}{2}.$$
  Thus, the general solution of \eqref{dform} is
  \begin{eqnarray*}
    D_{n+1} &=& k_{1}\lambda_{-}^{n+1}+k_{2}\lambda_{+}^{n+1} \\
     &=& k_{1}\left(\frac{B-\sqrt{B^2-4\tau_{n,n-1}\tau_{n-1,n}}}{2}\right)^{n+1}+k_{2}\left(\frac{B+\sqrt{B^2-4\tau_{n,n-1}\tau_{n-1,n}}}{2}\right)^{n+1}.
  \end{eqnarray*}
  When $n=0, D_{1}=0=B$. For $n=1$
  \begin{eqnarray*}
    D_{2}&=& k_{1}\left(\frac{B-\sqrt{B^2-4\tau_{0,1}\tau_{1,0}}}{2}\right)^{2}+k_{1}\left(\frac{B+\sqrt{B^2-4\tau_{0,1}\tau_{1,0}}}{2}\right)^{2}\\
    -\tau_{0,1}\tau_{1,0}&=& \frac{B\sqrt{B^2-\tau_{1,0}\tau_{0,1}}}{4}(k_{2}-k_{1})+\tau_{1,0}\tau_{0,1}(k_1+k_2).
  \end{eqnarray*}
  Comparing both sides of the last equation one obtains $k_2-k_1=0; k_2+k_1=-1$ which yields $k_{1}=k_{2}=-\frac{1}{2}.$
  Hence,
  \begin{equation}\label{dform3}
    D_{n+1}=-\frac{1}{2}\bigg[\left(\frac{B-\sqrt{B^2-4\tau_{n,n-1}\tau_{n-1,n}}}{2}\right)^{n+1}+\left(\frac{B+\sqrt{B^2-4\tau_{n,n-1}\tau_{n-1,n}}}{2}\right)^{n+1}\bigg].
  \end{equation}
  Applying \eqref{dform3} to $D_{n+1}=0$ yields $n+1$ values of the accessory parameter $B$.
  \end{proof}
  \section{Exact Solvability}\label{ESS}
  In this section, we establish exact solvability of the BHO.  By Turbiner (cf: Turbiner~\cite{TUB}, Lemma 2.2, p.13) a QES operator $T\in U_{s\ell(2,\mathbb{R})}$ has no terms $\mathcal{J}_{+}$ in positive grading if and only if it is an exactly solvable operator. Therefore setting $j=\frac{1}{2},$ the term $\mathcal{J}_{+}$ vanishes and $-H_{1}$ becomes
\begin{equation}\label{Hamj}
  -H_{e}=2\left(\mathcal{J}_{+}\mathcal{J}_{0}+\mathcal{J}_{0}\mathcal{J}_{+}\right)-\frac{g_2}{2}\left(\mathcal{J}_{0}\mathcal{J}_{-}+\mathcal{J}_{-}\mathcal{J}_{0}\right)-g_3(\mathcal{J}_{-})^2 +\frac{g_{2}}{4}\mathcal{J}_{-}+B.
\end{equation}
Explicitly, the eigenvalue differential equation associated with $H_{e}$ is given as
\begin{equation}\label{Hamexp}
  -H_{e}R(r):=\left((4r^3-g_2r-g_3)\frac{\vd^2 }{\vd r^2}+\frac{g_2}{4}\frac{\vd}{\vd r}+B\right)R(r)=0.
\end{equation}
The eigenfuntion $R_{m}(r)$ will be obtained using a combination of suitable gauge transformation and the technique of point canonical transformation (PCT) introduced by Shifman~\cite{SHIF} and  Levai~\cite{LEV} respectively. Gauge transformation allows the term in positive grading to be gauged away to obtain the gauge hamiltonian which may be solved exactly. PCT in its own right, allows the exactly solvability of the Hamiltonian in terms of a known orthogonal polynomial. The theorem stated below examines the exact solvability of $-H_{e}.$
\teo The radial wave function $R(r)$ which satisfies $-H_{e}R(r)=0$ is given using gauge transformation by
\[R(r)\approx k_{\pm}^{-\frac{1}{2}}\displaystyle{\prod_{s=1}^{3}}(r-e_s)^{-\nu_s}w_{\pm}(r)^{-\frac{\gamma\pm1}{2}}(w_{\pm}(r)-1)^{\frac{\nu-\gamma\pm1}{2}}P_{m}^{(\nu - \gamma,\gamma - 1)}(2w_{\pm}(r)-1),\]
where, $P_{m}^{\alpha,\beta}(\cdot)$ is a Jacobi polynomial and \[w_{\pm}(r)\sim \exp\left(\frac{g_2\pm\sqrt{g_2^2+64B(g_3+2)}}{4(g_3+2)}r\right).\]

\begin{proof}
To this end, let $R(r)$ be given as
\begin{equation}\label{pct}
  R(r)=f(r)F(w(r))
\end{equation}
 and $H_{e}$ be re-written in the form
\begin{equation}\label{gceq}
 \mathds{H}'=-\,\frac{1}{2}\frac{\vd^{2}}{\vd r^{2}}+\mathscr{A}(r)\frac{\vd}{\vd r}+\Delta V(r),
\end{equation}
 such that $V(r)\equiv \Delta V(r)-1/2 \mathscr{A}'(r)+1/2(\mathscr{A}(r))^{2} $ is a new exactly solvable potential obtained from $\mathds{H}'$. $\mathscr{A}(r)$ is called the \emph{gauge potential} and $\Delta V(r)$ is the modified potential from $\mathds{H}$ . The function $f(r)$ plays the role of a ground state eigenfunction of operator $\mathds{H}'$ in ~\eqref{gceq} and $F(w(r))$ is a special function satisfying a known second order differential equation
  \begin{equation}\label{sode}
    P(w)\ddot{F}(w)+Q(w)\dot{F}(w)+R(w)F(w)=0, \;\;\;\dot{}= \vd/\vd w
  \end{equation}
 of an orthogonal polynomial, where $P(w), Q(w)$ and $R(w)$ are polynomial functions in $w$. The gauge potential $\mathscr{A}(r)$ is determined  as
\begin{equation}\label{gpot}
  \mathscr{A}(r)=\frac{g_2}{8(4r^3-g_2r-g_3)}\equiv \frac{\nu_1}{r-e_1}+\frac{\nu_2}{r-e_2}+\frac{\nu_3}{r-e_{3}},
\end{equation}
where
\begin{eqnarray*}
  \nu_{1} &=& \frac{g_2}{32(e_1-e_2)(e_1-e_3)}; \\
  \nu_{2} &=& \frac{g_2}{32(e_2-e_1)(e_2-e_3)}; \\
   \nu_{3} &=& \frac{g_2}{32(e_3-e_1)(e_3-e_2)};
\end{eqnarray*}
and $e_1+e_2+e_3=0, g_2=4(e_1e_2+e_2e_3+e_1e_3), g_3=4e_1e_2e_3.$ The gauge potential given in \eqref{gpot} is used to obtain the phase function $$\displaystyle \Phi(r)=\int^{r} \mathscr{A}(u)\vd u=\ln \prod_{s=1}^{3}(r-e_s)^{\nu_{s}}.$$ Suppose $F(w)\;(\mathrm{where}\;w\equiv w(r))$ is in the class of orthogonal hypergeometric polynomials (see Bajpai~\cite{BAJ}) which satisfy the differential equation
\begin{equation}\label{GSHE}
  \mathds{\mathfrak{H}} F(w)=  w(w - 1)\ddot{F}(w)+[(\nu +1)w-\gamma]\dot{F}(w)-m(m + \nu)F(w)=0,
 \end{equation}
 where $\dot{}:=\frac{\vd}{\vd w}.$ Here and hereafter $m,\nu, \gamma\in\mathbb{Z}^{+}\cup \{0\}$ so that $F(w)$ may be written in terms of Jacobi polynomial of order $m$ as
 \begin{eqnarray*}
   F(w):=F_{m}(w) &=& _2F_1(-m,m+\nu; \gamma|w) \\
  &=&(-1)^{m}\frac{\Gamma(2m+\nu)m!}{\Gamma(m+\gamma)}P_{m}^{(\nu-\gamma,\gamma-1)}(2w-1).
 \end{eqnarray*}
The Jacobi polynomial $ P_{m}^{(\nu - \gamma, \gamma - 1)} (\cdot)$ (cf: Abramowitz and Stegun ~\cite{AS},\S 22.5-8, pp. 782-789) is given by
\begin{eqnarray*}
 P_{m}^{(\nu - \gamma,\gamma - 1)}(2w-1)&=&\left(
                              \begin{array}{c}
                                m+\nu - \gamma \\
                                 m\\
                              \end{array}
                            \right)F\left(-m, m+\nu,\nu - \gamma+1\bigg|1-w\right)  \\
                         &=& \frac{\Gamma(\nu - \gamma+m+1)}{m!\Gamma(\nu +m)}\sum_{n=0}^{m}\left(
                                                                                                  \begin{array}{c}
                                                                                                    m \\
                                                                                                    n \\
                                                                                                  \end{array}
                                                                                                \right)\frac{\Gamma(\nu+n+m)}{\Gamma(\nu - \gamma+n+1)}\left(1-w\right)^n.
\end{eqnarray*}
Assume that the operator $\mathbb{H}'$ is mapped into the hypergeometric operator $\mathds{\mathfrak{H}}$, that is, $\mathbb{H}'\mapsto\mathbb{\mathfrak{H}}$ by the variable transformation $R(r)=f(r)F(w(r))$ (using the point canonical transformation (PCT)) then
 the  function $f(r)$ is given by
\begin{eqnarray}\label{gst}
   f(r)&=&\frac{1}{\sqrt{w'(r)}}\exp\left(-\frac{1}{2}\int^{w(r)}\frac{Q(u)}{P(u)}\vd
   u\right)\exp\left(-\frac{1}{2}\Phi(r)\right)\nonumber\\
   &=& \frac{1}{\sqrt{w'(r)}}\exp\left(-\frac{1}{2}\int^{w(r)}\left(\frac{\gamma}{u}-\frac{\alpha + \beta -\gamma+ 1}{u-1}\right)\vd
   u\right)\exp\left(-\frac{1}{2}\Phi(r)\right)\nonumber\\
   &=&\exp\left(-\frac{1}{2}\Phi(r)\right)\sqrt{\frac{w(r)^{-\gamma}(w(r)-1)^{\alpha + \beta -\gamma+ 1}}{w'(r)}}\nonumber\\
   &=& \prod_{s=1}^{3}(r-e_s)^{\frac{\nu_{s}}{2}}w(r)^{-\frac{\gamma}{2}}(w(r)-1)^{\frac{\alpha + \beta -\gamma+ 1}{2}}(w'(r))^{-\frac{1}{2}}.\nonumber\\
\end{eqnarray}
The gauge function $f(r)$ satisfies the equation
\begin{equation}\label{pfr}
  \left\{P_{3}(r)\frac{\vd^{2}}{\vd r^{2}}+P_{2}(r)\frac{\vd}{\vd r}+ \Delta V_{num}-B\right\}f(r)=0
\end{equation}
where $\Delta V_{num}$ is the numerator of the potential $\Delta V.$ To evaluate the new variable $w(r)$, the expression needed is
\begin{eqnarray}\label{vej}
  \Delta V_{num}(r)-B &=& -P_2(r)\frac{f'(r)}{f(r)}-P_3(r)\frac{f''(r)}{f(r)} \nonumber\\
   &=&-P_{2}(r)\left\{\frac{1}{2}\frac{Q(w)w'(r)}{P(w)}-\frac{1}{2}\frac{w''(r)}{w'(r)}+\mathscr{A}(r)\right\}\nonumber\\
   && -P_{3}(r)\left\{\left(\frac{1}{2}\frac{Q(w)w'(r)}{P(w)}-\frac{1}{2}\frac{w''(r)}{w'(r)}+\mathscr{A}(r)\right)^{2}+\frac{1}{2}\frac{\dot{Q}(w)}{P(w)^{2}}\right.\nonumber\\
   &&\left.-\frac{1}{2}\frac{Q(w)\dot{P}(w)}{P(w)^{2}}-\frac{1}{2}\frac{w'''(r)}{w'(r)}-\frac{1}{2}\left(\frac{w''(r)}{w'(r)}\right)^{2}\right.\nonumber\\
   &&\left.-2\left(\frac{P'_{2}(r)}{P_{3}(r)}+\frac{P'_{3}(r)}{P_{3}(r)}\right)\mathscr{A}(r)\right\}.
\end{eqnarray}
To account for the constant part on the left hand side of equation~\eqref{vej}, an expression in terms of $w(r)$ must be equated to a constant on the right hand side. Let $\frac{w'}{w}=\mathrm{const.}$ (say $k_1<0$) then $w(r)=e^{k_1r}$ and the expression~\eqref{vej} yields
\begin{equation}\label{cst}
  -B \approx -\frac{g_2}{8}k_1+\frac{g_3}{4}k_1^2-\frac{g_3}{2}k_1^2-\frac{1}{2}k_1^2+ O(w(r))
\end{equation}
Equation~\eqref{cst} is re-written as a quadratic equation in $k_1$ as
\begin{equation}\label{qq1}
2(g_3+2)k_1^2+g_2k_1-8B\approx 0.
\end{equation}
and by general quadratic formula
\begin{equation}\label{qq2}
  k_1\equiv k_{\pm}\approx\;\frac{-g_2\pm\sqrt{g_2^2+64B(g_3+2)}}{4(g_3+2)}.
\end{equation}
Thus, $w(r)\equiv w_{\pm}(r)$ is given by \[w_{\pm}(r)\sim \exp\left(\frac{g_2\pm\sqrt{g_2^2+64B(g_3+2)}}{4(g_3+2)}r\right),\]
\begin{equation}\label{gsln}
f(r)\approx f_\pm(r)= k_{\pm}^{-\frac{1}{2}}\displaystyle{\prod_{s=1}^{3}}(r-e_s)^{-\frac{\nu_s}{2}}w_{\pm}(r)^{-\frac{\gamma\pm1}{2}}(w_{\pm}(r)-1)^{\frac{\nu-\gamma\pm1}{2}}
\end{equation}
  and
\begin{equation}\label{sol1}
R(r) \approx
                k_{\pm}^{-\frac{1}{2}}\displaystyle{\prod_{s=1}^{3}}(r-e_s)^{-\frac{\nu_s}{2}}w_{\pm}(r)^{-\frac{\gamma\pm1}{2}}(w_{\pm}(r)-1)^{\frac{\nu-\gamma\pm1}{2}}P_{m}^{(\nu - \gamma,\gamma - 1)}(2w_{\pm}(r)-1).
\end{equation}
We know (cf: Rusev~\cite{RUS}, \S I.3.1, p.12) that $P_{m}^{(\nu - \gamma,\gamma - 1)}(\cdot)$ is an orthogonal function with normalisation constant
$$\mathcal{N}_{m}=\frac{2^{\nu-1}\Gamma(m+\nu - \gamma+1)\Gamma(m+\gamma)}{m!(2m+\nu)\Gamma(m+\nu)}.$$
\end{proof}
\rem
The solutions expressed in \eqref{sol1} reveal that the values of $R_{m}(r)$ is a good approximation for the radial part of the eigenfunctions of BHE for sufficiently large values of $r$.
\section{Distributional Solution}\label{DSS}
In what follows, Fourier transform approach is employed to obtain the distributional solution of the radial part of BHO. Let $\Omega$ denote an open subset of $\mathbb{RP}^1\setminus\{e_1,e_2,e_3;\infty\},$ where $\mathbb{RP}^1$ is the real projective space and
$$\mathbb{RP}^1\setminus\{e_1,e_2,e_3;\infty\}=(-\infty,e_1)\cup(e_1,e_2)\cup (e_2,e_3)\cup(e_3,\infty)\subset \mathbb{R}.$$
Also, let $\mathscr{C}_{c}^{\infty}(\Omega)$ denote the space of infinitely differentiable  functions on $\Omega$ with compact supports. We assume that  $\mathscr{C}_{c}^{\infty}(\Omega)$ carries its usual inductive limit topology (cf: \cite{JBN}, \S2.2: 40).
	
Let $\chi_{\sigma}(r)=e^{i\sigma r}\;(\sigma\in\mathbb{R})$ be the characters of the group $\mathbb{R}$ and let $\chi_{-\sigma}(r)=\overline{\chi_{\sigma}(r)}=e^{-i\sigma r}$ be the complex conjugate of $\chi_{\sigma}$. The Fourier transform of any function $f\in L^{2}(\Omega,\omega(r)\vd r)$ is defined by the inner product $\langle\cdot,\cdot\rangle_\omega$ as
\begin{equation}\label{FFFT}
\langle f,\chi_{\sigma}\rangle_\omega=\int_{\Omega}f(r)\chi_{-\sigma}(r)\omega(r)\vd r,
\end{equation}
where $\vd\mu_\omega(r)=\omega(r)\vd r$ is a Radon measure on $\Omega$.

Let $\mathscr{D}'(\Omega)$ be the dual space of $\mathscr{C}_{c}^{\infty}(\Omega)$ endowed with strong dual topology. We note that
\begin{equation}\label{gelt}
\mathscr{C}_{c}^{\infty}(\Omega)\subset L^{2}(\Omega,\omega(r)\vd r)\subset\mathscr{D}'(\Omega)
\end{equation}
and the relative topology which the space $\mathscr{C}_c^{\infty}(\Omega)$ inherits from the Hilbert space topology defined by the norm $\|\cdot\|_{L^2}$ is strictly weaker than its usual inductive limit topology. In fact, $\mathscr{C}_c^{\infty}(\Omega)$ is dense in $L^{2}(\Omega,\omega(r)\vd r)$ with respect to the inherited Hilbert space topology.
It follows from \eqref{gelt}, that the Fourier transform of any function $\varphi(r)\in\mathscr{C}_{c}^{\infty}(\Omega),$ may be written as
\begin{equation}\label{FT1A}
\langle \varphi,\chi_{\sigma}\rangle_\omega=\int_{\Omega}\varphi(r)\chi_{-\sigma}(r)\omega(r)\vd r=\int_{\Omega}\varphi(r)e^{-i\sigma r}\omega(r)\vd r=\widehat{\varphi(r)}\equiv\mathscr{F}[\varphi(r)].
\end{equation}
Lemma~\ref{lemft} below is gleaned from Gel'fand and Shilov (\cite{GS}, Chapter II, \S2.2, pp.166-168).
\lem\label{lemft} Let $\delta(r)$ be a Dirac delta function and $P(r)$ be a polynomial of degree $\deg(P(r))=n$, then
\begin{description}
	\item[F1:] $\mathscr{F}[\delta^{(2m)}(r)]=(-i\sigma)^{2m}=(-1)^{m}\sigma^{2m};$
	\item[F2:]$\mathscr{F}[\delta^{(2m+1)}(r)]=(-i\sigma)^{2m+1}=(-1)^{m+1}\;i\sigma^{2m+1};$
	\item[F3:]$\mathscr{F}[P(r)]= \mathscr{F}[P(r)\cdot 1]=P\left(-i\frac{\vd}{\vd \sigma}\right)\mathscr{F}(1)= 2\pi P\left(-i\frac{\vd}{\vd \sigma}\right)\delta(\sigma);$ and
	\item[F4:]$	\mathscr{F}[P\left(\frac{\vd}{\vd r}\right)\delta(r)]= P(-i\sigma)\widehat{\delta(r)}=P(-i\sigma)\cdot 1=P(-i\sigma).$
\end{description}
Let
\begin{equation}\label{radop}
H^{g_2,g_3,n}=\left(4r^{3}- g_{2}r- g_{3}\right)\frac{\vd^{2}}{\vd r^{2}} -(n-\frac{1}{2})(6r^{2}-\frac{1}{2}g_2)\frac{\vd}{\vd r}+n(2n-1)r-B.
\end{equation}
\teo Consider the lemniscate case $g_2=1, g_3=0$ where $n$ are non-positive even numbers which describe radial BHO
  $$H^{g_2,g_3,n}=\left(4r^{3}- g_{2}r- g_{3}\right)\frac{\vd^{2}}{\vd r^{2}}-(n-\frac{1}{2})(6r^{2}-\frac{1}{2}g_2)\frac{\vd}{\vd r}+n(2n-1)r-B.$$
Its associated distributional solution is given by $R(r)=\sum_{k=1}^{\infty}a_{k}\delta^{(k)}(r)$ in $$Dom( H^{g_2,g_3,n})=\mathscr{C}_{c}^{\infty}(\Omega)\subset L^{2}(\Omega,\vd\mu_{\omega})(r)\subset\mathscr{D}'(\Omega)$$
 for $n=-2s$ provided that
\[a_{k}=K_{1}\left[\frac{\varsigma_{k,m}+\sqrt{\varsigma_{k,m}^2+4\epsilon_{k,m}}}{2}\right]^{k}+K_{2}\left[\frac{\varsigma_{k,m}-\sqrt{\varsigma_{k,m}^2+4\epsilon_{k,m}}}{2}\right]^{k},\]
\begin{eqnarray*}
 \varsigma_{k,m}&=&\frac{(1-2n)\big[3(k)_{m+1} +\frac{1}{4}(k)_{m-1}\big]}{\big[n(2n-1)(k)_{m+1} -(k)_{m}\cdot q\big]}\; \mathrm{and}\\
  \epsilon_{k,m} &=& \frac{\big[4(k)_{m+1} +(k)_{m-1}\big]}{\big[n(2n-1)(k)_{m+1} -(k)_{m}\cdot q\big]},
\end{eqnarray*}
where $K_{1},K_{2}\in(0,1)$ are arbitrary constants such that $K_{1}+K_{2}=1$ and hence
\begin{multline*}
 R(r)=\sum_{m=0}^{\lfloor~N_{1}\rfloor}\left[\delta(r)+\frac{(\varsigma_{1,m}+\Lambda_{1,m})}{2}\delta'(r)+\left(\frac{\epsilon_{2,m}(\varsigma_{1,m}+\Lambda_{1,m})}{2}+\varsigma_{2,m}\right)\delta''(r)\right.\\
 \left.+\left(\frac{(\epsilon_{3,m}\epsilon_{2,m}+\varsigma_{3,m})(\varsigma_{1,m}+\Lambda_{1,m})}{2}+\epsilon_{3,m}\varsigma_{2,m}\right)\delta'''(r)+\cdots\right].
\end{multline*}
with $\displaystyle\Lambda_{1,m}=(1-2K_{2})\sqrt{\varsigma_{1,m}^2+4\epsilon_{1,m}}.$
\begin{proof}
The domain $Dom(H)=\mathscr{C}_{c}^{\infty}(\Omega)$. Since the coefficients of the operator $H$ are smooth functions, then the space $\mathscr{C}_{c}^{\infty}(\Omega)$ is invariant under $H$. Let $\psi(r)\in\mathscr{C}_{c}^{\infty}(\Omega).$  As $\Omega\subset\overline{\Omega}$, then $\Omega$ is bounded and one has
 \begin{eqnarray*}
 \langle H\psi(r),\varphi\rangle_\omega&=& \int_{\Omega}\psi(r)(H\varphi(r))\vd\mu_\omega(r)\\
 &=&\langle\psi(r), H\varphi(r)\rangle_{\omega},\;\;\forall\; \varphi(r)\in\mathscr{C}_{c}^{\infty}(\Omega).
 \end{eqnarray*}
It follows that $H$ is a densely defined self-adjoint differential operator on the Hilbert space $L^{2}(\Omega,\vd\mu_{\omega}(r))$. Hence, $H$ is closed. This allows for extension of the domain of $H$ to cover the whole of $\mathscr{D}'(\Omega)$. It is clear that the range $Ran(H)=\mathscr{D}'(\Omega).$ Let the distributional solution of the differential equation defined by $H$, that is, equation~\eqref{rdq} be given as $\displaystyle R(r)=\sum_{k=0}^{\infty} a_k\delta^{(k)}(r)$, where $\displaystyle\;^{(k)}:=\frac{\vd^{k}}{\vd r^{k}}$. It is our task to determine the coefficients $a_{k}$. It follows from equation~\eqref{rdq} that
 \begin{equation}\label{ddf2}
   \langle HR(r), \chi_{\sigma} (r)\rangle_\omega=0.
 \end{equation}
Here and hereafter, let $\lfloor x\rfloor$ denote the monotone increasing \emph{floor function} which is the greatest integer greater or equal to $x\in\mathbb{R}^{+}$ (cf:\cite{KAM}, \S Eq.(36) :383) and let the weight function $\omega(r)$ be given as
 \begin{equation}\label{ddf3}
   \omega(r)=(r-e_1)^{N_1}(r-e_2)^{ N_2}(r-e_3)^{N_3},
 \end{equation}
 We know (Abramowitz and Stegun~1972,\S 18.14: 658)  that the values of the roots of Weierstrass polynomial $4r^{3}-g_{2}r-g_{3}$ for the Lemniscatic case $(g_2=1, g_3=0)$ at half-periods are $e_1=0,\; e_2=\frac{1}{2}$ and $\;e_3=-\frac{1}{2}.$  By substituting these values into $N_i$ ($i=1,2,3$) gives
  \begin{eqnarray*}
   N_1 &=& -\frac{(n-\frac{1}{2})(6e_1^2-\frac{g_2}{2})}{4(e_1-e_2)(e_1-e_3)}-1=-\frac{2n+3}{4}, \\
   N_2 &=&- \frac{(n-\frac{1}{2})(6e_2^2-\frac{g_2}{2})}{4(e_2-e_1)(e_2-e_3)}-1=-\frac{2n+3}{4} ,\\
   N_3 &=& -\frac{(n-\frac{1}{2})(6e_3^2-\frac{g_2}{2})}{4(e_3-e_2)(e_3-e_1)}-1=-\frac{2n+3}{4}.
 \end{eqnarray*}
 Since, $N_1=N_2=N_3=-\frac{2n+3}{4},$ the function $\omega(r)$ in equation~\eqref{ddf3} may be expressed as
 \begin{equation}\label{ddf4}
   \omega(r)=(r^3-\frac{1}{4}r)^{N_1}=\sum_{p=0}^{N_1}\binom{N_1}{p}r^{3(N_1-p)}\left(-\frac{1}{4}\right)^{p}r^p=\sum_{p=0}^{N_1}(-1)^p2^{-2p}\binom{N_1}{p}r^{3N_1-2p}
 \end{equation}
 Following van der Waall~(\cite{vandw}, Chapter 4, p.55), $n\in\mathbb{Z}$ can assume a negative value since it is obtained from the Lam\'{e} operator with reducible monodromy. We deduce that since $N_{i}(i=1,2,3)$ depends on $n$ we can choose $n$ such that $n=-2s,s=1,2,\ldots\;.$ This implies that $N_{i}=-\frac{-2(2s)+3}{4}=\frac{4s-3}{4}=s-\frac{3}{4}.$

Let $g_2=1,g_3=0,n=-2s, s=1,2,\ldots$. Then, the Hamiltonian in this case is
 \begin{equation}\label{fob13}
   H_{1}=H^{g_2=1,g_3=0,n}:=(4r^3-r)\frac{\vd^2}{\vd r^2}-(n-\frac{1}{2})(6r^2-\frac{1}{2})\frac{\vd}{\vd r}+n(2n-1)r-B.
 \end{equation}
  Therefore, by Fourier transform
 \begin{equation}\label{fob14}
   \langle H_1 R(r),\chi_{\sigma}(r)\rangle_{\omega}=0.
 \end{equation}
 This can further be viewed as
 \begin{equation}\label{fob15}
 \sum_{p=0}^{N_1}(-1)^p2^{-2p}\binom{N_1}{p}\bigg\langle r^{3N_1-2p} H_1R(r),\chi_{\sigma}(r)\bigg\rangle=0.
 \end{equation}
 Explicitly,  equation~\eqref{fob15} becomes
 \begin{multline}\label{fob16}
    \sum_{p=0}^{N_1}(-1)^p2^{-2p}\binom{N_1}{p}\\
    \times\bigg\langle r^{3N_1-2p}\left[(4r^3-r)\frac{\vd^2}{\vd r^2}-(n-\frac{1}{2})(6r^2-\frac{1}{2})\frac{\vd}{\vd r}+n(2n-1)r-B\right]R(r),\chi_{\sigma}(r)\bigg\rangle=0.
 \end{multline}
 Considering each independent term and ignoring the common multiplier $2\pi$ because of the polynomial coefficients of the differential operator which easily factors out, we get
 \begin{eqnarray}
 \sum_{p=0}^{N_1}&(-1)^p&2^{-2p}\binom{N_1}{p}\mathscr{F}[(4r^{3+3N_1-2p}-r^{1+3N_1-2p})R''(r)]\nonumber\\
 &=& \sum_{p=0}^{N_1}(-1)^p2^{-2p}\binom{N_1}{p} \langle (4r^{3+3N_1-2p}-r^{1+3N_1-2p})R''(r),\chi_{\sigma}(r)\rangle\nonumber\\
    &=& \sum_{p=0}^{N_1}\sum_{k=0}^{\infty}(-1)^p2^{-2p}\binom{N_1}{p}a_{k}\nonumber\\
    && \cdot\bigg\langle (4 r^{3+3N_1-2p}-r^{1+3N_1-2p})\delta^{(k+2)}(r),\chi_{\sigma}(r)\bigg\rangle \nonumber\\
   &=& \sum_{p=0}^{N_1}\sum_{k=0}^{\infty}(-1)^p2^{-2p}\binom{N_1}{p}a_{k}\nonumber\\
    && \cdot \left(4(-i)^{3+3N_1-2p}\left(\frac{\vd}{\vd\sigma}\right)^{3+3N_1-2p}-(-i)^{1+3N_1-2p}\left(\frac{\vd}{\vd\sigma}\right)^{1+3N_1-2p}\right)(-i)^{k+2}\sigma^{k+2} \nonumber\\
    &=& \sum_{p=0}^{N_1}\sum_{k=0}^{\infty}(-1)^p2^{-2p}\binom{N_1}{p}a_{k}(-i)^{3N_1-2p+k}\nonumber\\
    && \cdot (-i)\left(4(k+2)_{3+3N_1-2p}+(k+2)_{1+3N_1-2p}\right);\label{fod1}
    \end{eqnarray}
 \begin{eqnarray}
-(n-\frac{1}{2})\sum_{p=0}^{N_1}&(-1)^p&2^{-2p}\binom{N_1}{p}\mathscr{F}[(6r^{2+3N_1-2p}-\frac{1}{2})r^{3N_1-2p}R'(r)] \nonumber\\
&=&-\sum_{p=0}^{N_1}(-1)^p2^{-2p}\binom{N_1}{p}(n-\frac{1}{2})\nonumber\\
&&\cdot\bigg\langle  \left(6r^{2+3N_1-2p}-\frac{1}{2}r^{3N_1-2p}\right)R'(r),\chi_{\sigma}(r)\bigg\rangle \nonumber\\
    &=&\sum_{p=0}^{N_1}\sum_{k=0}^{\infty}(-1)^p2^{-2p}\binom{N_{1}}{p}\nonumber\\
    &&\cdot a_{k}\bigg\langle \left(3(1-2n) r^{2+3N_1-2p}+\frac{1}{2}(n-\frac{1}{2})r^{3N_1-2p}\right)\delta^{(k+1)}(r),\chi_{\sigma}(r)\bigg\rangle\nonumber\\  &=&\sum_{p=0}^{N_1}\sum_{k=0}^{\infty}(-1)^p2^{-2p}\binom{N_{1}}{p}(1-2n) a_{k}
\left(3 (-i)^{2+3N_1-2p} \left(\frac{\vd}{\vd\sigma}\right)^{2+3N_1-2p}\right.\nonumber\\
   && \left.-\frac{1}{4}(-i)^{3N_1-2p}\left(\frac{\vd}{\vd\sigma}\right)^{3N_1-2p}\right)(-i)^{k+1}\sigma^{k+1} \nonumber\\
      &=&\sum_{p=0}^{N_1}\sum_{k=0}^{\infty}(-1)^p2^{-2p}\binom{N_{1}}{p}(1-2n)(-i)^{3N_1-2p+k}\nonumber\\
    &&\cdot a_{k}i \left(3(k+1)_{2+3N_1-2p}+\frac{1}{4}(k+1)_{3N_1-2p}\right)   ; \label{fod2}
    \end{eqnarray}
\begin{eqnarray*}
  \sum_{p=0}^{N_1}&(-1)^p&2^{-2p}\binom{N_1}{p}\mathscr{F}[\left(n(2n-1)r^{1+3N_{1}}-Br^{3N_{1}-2p}\right)R(r)] \nonumber\\
  &=&\sum_{p=0}^{N_1}(-1)^p2^{-2p}\binom{N_1}{p}  \left\langle\left(n(2n-1)r^{1+3N_{1}-2p}-Br^{3N_{1}-2p}\right)R(r),\chi_{\sigma}(r)\right\rangle\nonumber\\
     &=& \sum_{p=0}^{N_1}(-1)^p2^{-2p}\binom{N_1}{p}  \sum_{k=0}^{\infty}a_{k}\nonumber\\
    &&\cdot\bigg\langle \left(n(2n-1) r^{1+3N_{1}-2p}-Br^{3N_{1}-2p}\right)\delta^{(k+1)}(r),\chi_{\sigma}(r)\bigg\rangle\nonumber\\
   &=& \sum_{p=0}^{N_1}(-1)^p2^{-2p}\binom{N_1}{p}  \sum_{k=0}^{\infty}a_{k}(-i)^{3N_{1}-2p}\nonumber\\
   &&\left(n(2n-1)(-i) \left(\frac{\vd}{\vd\sigma}\right)^{1+3N_{1}-2p}-B\left(\frac{\vd}{\vd\sigma}\right)^{3N_{1}-2p}\right)(-i)^{k}\sigma^{k}
    \end{eqnarray*}
    \begin{eqnarray}
   &=& \sum_{p=0}^{N_1}(-1)^p2^{-2p}\binom{N_1}{p}  \sum_{k=0}^{\infty}a_{k}(-i)^{3N_{1}-2p+k}\left(n(2n-1)(-i) (k)_{1+3N_{1}-2p}-B(k)_{3N_{1}-2p}\right).\nonumber\\ \label{fod3}
   \end{eqnarray}
Re-coupling equations~\eqref{fod1}-\eqref{fod3}, $\mathscr{F}_{\omega}[H_1R(r)]=0$ explicitly becomes
 \begin{multline}\label{fod7}
\sum_{p=0}^{N_1}\sum_{k=0}^{\infty}(-1)^p2^{-2p}\binom{N_1}{p}(-i)^{3N_1-2p+k}a_{k}\bigg\{4i(k+2)_{3+3N_1-2p}+i(k+2)_{1+3N_1-2p}\\
 +(1-2n)i\left(3(k+1)_{2+3N_1-2p}+\frac{1}{4}(k+1)_{3N_1-2p}\right)\\
+ \left(n(2n-1)(-i) (k)_{1+3N_{1}-2p}-B(k)_{3N_{1}-2p}\right)\bigg\}=0.
 \end{multline}
Multiplying through equation~\eqref{fod7} by $(-1)^{p}(i)^{3N_1-2p+k}$ one obtains
\begin{multline}\label{fod8}
\sum_{p=0}^{N_1}\sum_{k=0}^{\infty}2^{-2p}\binom{N_1}{p}a_{k}\bigg\{4i(k+2)_{3+3N_1-2p}+i(k+2)_{1+3N_1-2p}\\
 +(1- 2n)i\left(3(k+1)_{2+3N_1-2p}+\frac{1}{4}(k+1)_{3N_1-2p}\right)\\
+ \left(n(2n-1)(-i) (k)_{1+3N_{1}-2p}-B(k)_{3N_{1}-2p}\right)\bigg\}=0.
\end{multline}
Multiplying through equation~\eqref{fod8} by $-i$ gives
\begin{multline}\label{fod9}
\sum_{p=0}^{N_1}\sum_{k=0}^{\infty}2^{-2p}\binom{N_1}{p}a_{k}\bigg\{4(k+2)_{3+3N_1-2p}+(k+2)_{1+3N_1-2p}\\
 -(2n-1)\left(3(k+1)_{2+3N_1-2p}+\frac{1}{4}(k+1)_{3N_1-2p}\right)\\
+ \left(-n(2n-1) (k)_{1+3N_{1}-2p}+iB(k)_{3N_{1}-2p}\right)\bigg\}=0.
\end{multline}
 In what follows, we will study the nature of solution in terms of some comb-like function so as to avoid fractional derivatives (see Poularikas~\cite{PAD}, \S7.21.1). When $k<3\lfloor~N_{1}\rfloor-2p,$ one obtains a trivial solution. On the other hand, for $k+2\geqslant 3\lfloor~N_{1}\rfloor-2p$ or $k\geqslant 3\lfloor~N_{1}\rfloor-2p-2.$ Let $m=3\lfloor~N_{1}\rfloor-2p$ then $k\geqslant m-2$ and when $p=0, m=3\lfloor~N_{1}\rfloor.$ This makes equation~\eqref{fod9} to take the form
\begin{multline}\label{fod10}
\sum_{p=0}^{\lfloor~N_1\rfloor}\sum_{k=0}^{\infty}2^{-2p}\binom{N_1}{p}a_{k}\bigg\{4(k+2)_{m+3}+(k+2)_{m+1}\\
 -(2n-1)\left(3(k+1)_{m+2}+\frac{1}{4}(k+1)_{m}\right)-\left(n(2n-1) (k)_{m+1}+iB(k)_{m}\right)\bigg\}=0.
\end{multline}
By shifting index and setting $q=-iB\in\mathbb{C},$ we get
\begin{multline}\label{fod11}
\sum_{p=0}^{\lfloor~N_1\rfloor}\sum_{k=0}^{\infty}2^{-2p}\binom{\lfloor~N_1\rfloor}{p}a_{k}\bigg\{4(k+2)_{m+3}+(k+2)_{m+1}\\
 -(2n-1)\left(3(k+1)_{m+ 2}+\frac{1}{4}(k+1)_{m}\right)-\left(n(2n-1) (k)_{m+1}+q(k)_{m}\right)\bigg\}=0.
\end{multline}
We remark here that if $B$ is purely imaginary then $q\in \mathbb{R}.$ If $B$ is a real number then $q=0.$ If $B\in\mathbb{C}$ with $\Im B\neq 0$ then the real part of $B,\Re B$ must be equal to zero and hence $q\in\mathbb{R}.$

From equation~\eqref{fod11}, one obtains the recurrence equation
\begin{multline*}
  a_{k-2}\big[4(k)_{m+1} +(k)_{m-1}\big]+ a_{k-1}\big[3(1-2n)(k)_{m+1} +\frac{1-2n}{4}(k)_{m-1}\big]\\
 -a_{k}\big[n(2n-1)(k)_{m+1} -(k)_{m}\cdot q\big]=0.
\end{multline*}
Thus, for $k\geqslant 2,n=-2s,s=1,2,\ldots$
\begin{equation}\label{fod12}
  a_{k}=\;\frac{\big[4(k)_{m+1} +(k)_{m-1}\big]a_{k-2}+(1-2n)\big[3(k)_{m+1} +\frac{1}{4}(k)_{m-1}\big]a_{k-1}}{\big[n(2n-1)(k)_{m+1} -(k)_{m}\cdot q\big]}.
\end{equation}
In what follows, the method of characteristic polynomial is used in obtaining solution of the recurrence equation (see Lipschutz and Lipson~\cite{LILI}, Chapter 6, \S 6.7, pp.113-114 ).
To solve the recurrence equation, we set $a_{k}=t^{k}$ so that $a_{0}=1.$ Also setting
\begin{eqnarray*}
  \varsigma_{k,m}&=&\frac{(1-2n)\big[3(k)_{m+1} +\frac{1}{4}(k)_{m-1}\big]}{\big[n(2n-1)(k)_{m+1} -(k)_{m}\cdot q\big]}\; \mathrm{and}\\
  \epsilon_{k,m} &=& \frac{\big[4(k)_{m+1} +(k)_{m-1}\big]}{\big[n(2n-1)(k)_{m+1} -(k)_{m}\cdot q\big]}
\end{eqnarray*}
we have
\[t^{k}-\varsigma_{k,m}t^{k-1}-\epsilon_{k,m}t^{k-2}=0.\]
Dividing through by $t^{k-2}$ we get
\[t^{2}-\varsigma_{k,m}t-\epsilon_{k,m}=0.\]
This is solved by quadratic formula to get the roots
\[t_{\pm}=\frac{\varsigma_{k,m}\pm\sqrt{\varsigma_{k,m}^2+4\epsilon_{k,m}}}{2}.\]
Thus, since the discriminant $D=\varsigma_{k,m}^2+4\epsilon_{k,m}>0$ we have two distinct real values given by
\[a_{k}=K_{1}\left[\frac{\varsigma_{k,m}+\sqrt{\varsigma_{k,m}^2+4\epsilon_{k,m}}}{2}\right]^{k}+K_{2}\left[\frac{\varsigma_{k,m}-\sqrt{\varsigma_{k,m}^2+4\epsilon_{k,m}}}{2}\right]^{k}.\]

Using the initial condition $k=0$ we get $a_{0}=K_{1}+K_{2}=1.$ When $k=1,$
\begin{eqnarray*}
  a_{1} &=& (K_{1}+K_{2})\frac{\varsigma_{1,m}}{2} +(K_{1}-K_{2})\frac{\sqrt{\varsigma_{1,m}^2+4\epsilon_{1,m}}}{2}\\
  a_{1} &=& \frac{\varsigma_{1,m}}{2}+(1-2K_{2})\frac{\sqrt{\varsigma_{1,m}^2+4\epsilon_{1,m}}}{2} \;\;\because K_{1}+K_{2}=1.
\end{eqnarray*}
Plugging $a_0=1,a_{1}=\frac{\varsigma_{1,m}}{2}+(1-2K_{2})\frac{\sqrt{\varsigma_{1,m}^2+4\epsilon_{1,m}}}{2}$ into equation~\eqref{fod12} we get
\begin{eqnarray*}
  a_{2} &=&\epsilon_{2,m}a_{1}+\varsigma_{2,m}a_{0} \\
   &=& \epsilon_{2,m} \left[\frac{\varsigma_{1,m}}{2}+(1-2K_{2})\frac{\sqrt{\varsigma_{1,m}^2+4\epsilon_{1,m}}}{2}\right]+\varsigma_{2,m}.
\end{eqnarray*}
Also, plugging in $a_0,a_1,a_2$ into recurrence equation $a_{k}=\epsilon_{k,m}a_{k-1}+\varsigma_{k,m}a_{k-2}$
gives
\begin{eqnarray*}
  a_{3} &=&\epsilon_{3,m}a_{2}+\varsigma_{3,m}a_{1}  \\
   &=&\epsilon_{3,m}\left[\epsilon_{2,m} \left(\frac{\varsigma_{1,m}}{2}+(1-2K_{2})\frac{\sqrt{\varsigma_{1,m}^2+4\epsilon_{1,m}}}{2}\right)+\varsigma_{2,m}\right]\\
   &&\;\;+\varsigma_{3,m}\left[\frac{\varsigma_{1,m}}{2}+(1-2K_{2})\frac{\sqrt{\varsigma_{1,m}^2+4\epsilon_{1,m}}}{2}\right]\\
   &=&\frac{\epsilon_{3,m}\epsilon_{2,m}\varsigma_{1,m}}{2}+\epsilon_{3,m}\varsigma_{2,m}+\frac{\varsigma_{3,m}\varsigma_{1,m}}{2}+\frac{(1-2K_{2})(\epsilon_{3,m}\epsilon_{2,m}+\varsigma_{3,m})\sqrt{\varsigma_{1,m}^2+4\epsilon_{1,m}}}{2}.
\end{eqnarray*}
Here, $K_{2}\in (0,1).$ By setting $\displaystyle\Lambda_{1,m}=(1-2K_{2})\sqrt{\varsigma_{1,m}^2+4\epsilon_{1,m}},$  the radial distributional solution is given by
\begin{multline*}
 R(r)=\sum_{m=0}^{\lfloor~N_{1}\rfloor}\left[\delta(r)+\frac{(\varsigma_{1,m}+\Lambda_{1,m})}{2}\delta'(r)+\left(\frac{\epsilon_{2,m}(\varsigma_{1,m}+\Lambda_{1,m})}{2}+\varsigma_{2,m}\right)\delta''(r)\right.\\
 \left.+\left(\frac{(\epsilon_{3,m}\epsilon_{2,m}+\varsigma_{3,m})(\varsigma_{1,m}+\Lambda_{1,m})}{2}+\epsilon_{3,m}\varsigma_{2,m}\right)\delta'''(r)+\cdots\right].
\end{multline*}
It is obvious here also that whenever $K_{2}=\frac{1}{2}, \Lambda_{1,m}=0, \;\forall\; m=0,\ldots, \lfloor N_{1}\rfloor.$
\end{proof}
\section{Conclusion}
In this paper, it is clearly seen that the radial part of a complex differential operator can be obtained using asymptotic variable transformation method. The BHO was written as an element in the center of universal enveloping algebra of the Lie group $SL(2,\mathbb{R}).$ The quasi-exact solvability and exact solvability of the radial part of BHO has been examined. The distributional solution of the radial part of BHO has been obtained. The three solutions is a clear description of the relationship of the Gel'fand triple.
\section*{Declarations}
\begin{itemize}
\item\textbf{Compliance with Ethical Standards}: The authors adhered to all ethical standards for the publication of this paper.
\item\textbf{Author's Contribution}: U.S. Idiong reviewed the literature, computed the results and carried out the typesetting of this paper. U.N. Bassey proposed the problem of evaluating radial distributions associated with differential operators that define Fuchsian equations. O.S. Obabiyi critiqued the results. 
  \item \textbf{Conflict of Interest}: The authors declare that there are no conflicts of interest in the publication of this paper.
  \item \textbf{Funding}: This research was partly funded by TETfund Nigeria.
 \item \textbf{Data Availability Statements}: Not applicable.
\item\textbf{Ethical Conduct}: This paper is not under consideration for publication in any other journal.
\item\textbf{Acknowledgement}: We acknowledge Fritz Gesztesy, Ram P. Kanwal and other authors whose works have served as a foundation on which these results have been built.
\end{itemize}

\end{document}